\documentclass[12pt,openany]{article} 

\usepackage[utf8]{inputenc} \usepackage{kpfonts}

\usepackage[spanish,es-noquoting]{babel} 
	\addto\captionsspanish{}	
	\addto\captionsspanish{}	
	\addto\captionsspanish{}

\usepackage[numbers]{natbib}
\usepackage{babelbib}

\usepackage{geometry}                         
\usepackage{amsmath}       
\usepackage{amsfonts}     	
\usepackage{amssymb}

\usepackage{bbold}

\usepackage{graphicx}     	

\usepackage[lofdepth,lotdepth]{subfig}

\usepackage{tikz}

\usepackage[american]{circuitikz}

\usepackage{siunitx}	

\usepackage[colorlinks=true,urlcolor=blue,linkcolor=blue,citecolor=red]{hyperref}
\usepackage{float}	

\usepackage{booktabs}		

\usepackage{multicol}

\usepackage{listings}
\usepackage[numbered,autolinebreaks]{mcode}

\usepackage{lastpage}

\usepackage{fancyhdr}
\usepackage{marvosym}

\usepackage{blindtext}

\usepackage[abs]{overpic}

\newcommand{\KP}[1]{%
	\begin{tikzpicture}[baseline=-\dimexpr\fontdimen22\textfont2\relax]
		#1
	\end{tikzpicture}%
}

\newcommand{\KPB}{%
	\KP{
		\draw[color=gray,thick] (-0.3,0.3) -- (0.3,-0.3);
		\draw[color=gray,thick] (-0.3,-0.3) -- (-0.05,-0.05);
		\draw[color=gray,thick] (0.05,0.05) -- (0.3,0.3);
	}%
}
\newcommand{\KPC}{%
	\KP{%
		\draw[color=gray,thick] (-0.3,0.3) .. controls (0,-0.05) .. (0.3,0.3);
		\draw[color=gray,thick] (-0.3,-0.3) .. controls (0,0.05) .. (0.3,-0.3);
	}%
}
\newcommand{\KPD}{%
	\KP{%
		\draw[color=gray,thick] (-0.3,-0.3) .. controls (0.05,0) .. (-0.3,0.3);
		\draw[color=gray,thick] (0.3,-0.3) .. controls (-0.05,0) .. (0.3,0.3);
	}%
}

\newcommand\quotient[2]{
	\mathchoice
	{
		\text{\raise1ex\hbox{$#1$}\Big/\lower1ex\hbox{$#2$}}%
	}
	{
		#1\,/\,#2
	}
	{
		#1\,/\,#2
	}
	{
		#1\,/\,#2
	}
}

\newtheorem{theorem}{Theorem}[section]

\newtheorem{definition}[theorem]{Definici\'on}

\newtheorem{example}[theorem]{Ejemplo}

\newtheorem{corollary}[theorem]{Corollary}

\begin{document}
\author{Gabriel Montoya-Vega}

\title{Una Mirada Inicial a la Teor\'ia de Nudos y a la Homolog\'ia de Khovanov}


\date{} 

\pdfbookmark[1]{Portada}{portada} 	
\maketitle							


\begin{abstract}
La teor\'ia matem\'atica de nudos estudia las incrustaciones de c\'irculos en el espacio $\mathbb{R}^3$, siendo uno de los problemas fundamentales la clasificaci\'on. La introducci\'on de teor\'ias de homolog\'ia produce estructuras matem\'aticas complejas generando nuevas oportunidades de investigaci\'on. Con motivo del \textit{Encuentro Internacional de Matem\'aticas (EIMAT)} a realizarse en la Universidad del Atl\'antico en Barranquilla, Colombia en Noviembre de 2023, en este art\'iculo a manera de exposici\'on brindamos una primera mirada a la homolog\'ia de Khovanov y a la sucesi\'on larga de Khovanov. Adem\'as se presenta un resumen de los or\'igenes hist\'oricos de la teor\'ia que nos puede llevar incluso al a\~no 2600 a.e.c., pasando por la Italia del siglo XV y la Escocia del siglo XIX y se dan referencias para lecturas m\'as detalladas. Adicional a discutir la construcci\'on de la homolog\'ia de Khovanov a partir del polinomio bracket, el objetivo principal de esta publicaci\'on es atacar la barrera del idioma (una traducci\'on al ingl\'es tambi\'en ser\'a publicada) y fomentar el estudio de la teor\'ia de nudos y la homolog\'ia de Khovanov en Colombia y Latinoam\'erica en general.
\vspace{1.5 cm}

\begin{center}
	\textbf{Abstract}\\
\end{center}

The mathematical theory of knots studies the embeddings of circles into the space $\mathbb{R}^3$, being the classification one of the fundamental problems. The introduction of homology theories results in complex mathematical structures that generate new research opportunities. On occasion of the \textit{Encuentro Internacional de Matem\'aticas (EIMAT) (International Meeting of Mathematics)} to be celebrated at the Universidad del Atl\'antico in Barranquilla, Colombia in November 2023, in this article, in an expository way, we offer a first look into Khovanov homology and the long exact sequence of Khovanov homology. Moreover, we present a summary of the historical origins of the theory which can take us as early as the year 2600 BCE, passing through Italy of the XV century, Scotland of the XIX century, and we give references for further and more detailed discussions. Additionally to showing the construction of Khovanov homology from the Kauffman bracket polynomial, the main objective in publishing this article is to attack the language barrier and popularize knot theory and Khovanov homology in Colombia and Latin-America in general.
\end{abstract}

\vspace{1 cm}
Palabras claves: \textit{Nudos y enlaces, teor\'ia de nudos, polinomio bracket, homolog\'ia de Khovanov, sucesi\'on larga de la homolog\'ia de Khovanov, nudos toroidales.}\\
\indent Clasificaci\'on AMS: 01A05, 57K10, 57K14, 57K18


\tableofcontents
\listoffigures
\listoftables

\newpage
\section{Introducci\'on}
La teor\'ia matem\'atica de nudos es un campo fascinante de investigaci\'on que envuelve distintas \'areas de la ciencia como la matem\'atica, f\'isica, qu\'imica, entre otras. La principal motivaci\'on para publicar este art\'iculo es el \textit{Encuentro Internacional de Matem\'aticas (EIMAT)} a realizarse en la Universidad del Atl\'antico en Barranquilla, Colombia en Noviembre de 2023. El art\'iculo est\'a organizado de la siguiente manera. En la Secci\'on \ref{HIST} se brinda un resumen de la historia temprana de la teor\'ia. En la Secci\'on \ref{INTROToN} presentamos las primeras definiciones formales de la teor\'ia, entre ellos los diagramas de nudos y los movimientos de Reidemeister. En la Secci\'on \ref{HOM} se construye la homolog\'ia de Khovanov partiendo del polinomio bracket y en la Secci\'on \ref{suclarga} se construye la sucesi\'on larga de la homolog\'ia de Khovanov. Finalmente en la Secci\'on \ref{HOMTORUS} se calcula la homolog\'ia de los nudos toroidales $T(2,n)$ y en la Secci\'on \ref{CONCLU} se presentan posibles futuras direcciones de investigaci\'on.

\section{Historia}\label{HIST}

La teor\'ia de nudos, hist\'oricamente considerada como sub\'area de la topolog\'ia, es actualmente un campo de investigaci\'on matem\'atico que se ha desarrollado intensamente hasta alcanzar una alta independencia. Adem\'as de tener una fuerte conexi\'on con otras \'areas como el \'algebra, la combinatoria y la mec\'anica estad\'istica, goza de una historia fascinante. En esta secci\'on presentamos un resumen de la historia de la teor\'ia matem\'atica de nudos y proveemos referencias para lecturas mas detalladas. Esta historia consiste de un viaje en el tiempo que nos puede llevar incluso hasta el a\~no 2600 a.e.c. pasando por la Italia del siglo XV y por la Escocia del siglo XIX. Existe evidencia clara de que la humanidad ha mostrado inter\'es en los nudos desde tiempos muy tempranos en la historia, de hecho podr\'iamos solo especular que llev\'o a comerciantes en Anatolia (Asia menor) en el a\~no 1700 a.e.c a usar nudos y trenzas en sus sellos y firmas. Para mas informaci\'on al respecto y gr\'aficas relacionadas, referirse a \cite{Col, Prz}. Escavaciones hechas en la regi\'on de Lerna por la Escuela Americana de Estudios Cl\'asicos, bajo la direcci\'on del Profesor J. L. Caskey, encontraron varias impresiones de sellos en arcilla, muchos de ellos que corresponden al a\~no 2200 a.e.c. contienen nudos y enlaces \cite{Hea}. Es posible encontrar evidencia a\'un mas temprana (2600-2500 a.e.c) de sellos conteniendo nudos en \cite{Wolk}. Damos un salto en el tiempo hasta el siglo I donde un f\'isico griego llamado Heraklas, escribi\'o un ensayo donde daba instrucciones paso a paso de como atar 18 cabestrillos ortop\'edicos. Este trabajo sobrevivi\'o gracias a que Oribasius de P\'ergamo (m\'edico del emperador Juliano el Ap\'ostata) lo incluy\'o a finales del siglo IV en el escrito \textit{Sinagogas m\'edicas}, el cual inclu\'ia todo el saber m\'edico de la \'epoca. El escrito mas reciente de \textit{Sinagogas m\'edicas} data del siglo X y fue  llevado a Italia por el griego J. Lascaris, un refugiado de Constantinopla. El escrito de Heraklas no conten\'ia ilustraciones y s\'olo hasta el a\~no 1500 un artista an\'onimo dibuj\'o uno de sus nudos en el manuscrito de Oribasius de P\'ergamo. Una traducci\'on al latin realizada por Vidus Vidius (m\'edico del rey Francisco de Francia I) contiene illustraciones hechas por el artista italiano Francesco Primaticcio; ver por ejemplo \cite{PBIMW}. Otros artistas del renacimiento realizaron trabajos que de alguna forma est\'an relacionados a los nudos como Albrecht D\"urer y Leonardo da Vinci. Es importante mencionar que la historia de los escritos de Heraklas son importantes en la teor\'ia de nudos aunque sean s\'olo aplicaciones.

\ 

Podemos decir que la topolog\'ia y la teor\'ia de nudos moderna tienen su origen gracias a Gottfried Wilhelm Leibniz. Hay evidencia de una carta de Leibniz enviada a Christiaan Huygens en 1679, en la que habla de la necesidad de un nuevo tipo de an\'alisis. Leibniz lo llam\'o \textit{Analysis situs} o \textit{Geometr\'ia de la posici\'on}. En la carta escribi\'o: 
\begin{center}    
\textit{No estoy satisfecho con el \'algebra, ya que no produce las demostraciones m\'as cortas ni las construcciones m\'as hermosas de la geometr\'ia. En consecuencia, en vista de esto, considero que necesitamos otro tipo de an\'alisis, geom\'etrico o lineal, que trate directamente con la posici\'on, como el \'algebra trata con la magnitud.}
\end{center}

De cualquier manera, no es completamente claro si Leibniz ten\'ia alg\'un ejemplo de un problema que perteneciera al \textit{Analysis situs}. De hecho, el primer ejemplo conocido es el famoso problema de los puentes de K\"onigsberg. Este problema fue propuesto por Heinrich K\"uhn en 1735. K\"uhn comunic\'o el problema a Leonhard Euler sugiriendo que ser\'ia relacionado a la geometr\'ia de la posici\'on. Aunque en primeras instancias podr\'iamos decir que Euler no estaba convencido de la utilidad o importancia de este problema, en 1736 resulta publicando un art\'iculo donde detalla su soluci\'on \cite{Eul}. El art\'iculo titulado \textit{Solutio problematis ad geometriam situs pertinentis} (Soluci\'on de problemas relacionados con la geometr\'ia de la posici\'on), da origen a la teor\'ia de grafos y a la topolog\'ia.  El nacimiento de la teor\'ia de nudos se da cuando el matem\'atico Franc\'es Alexandre-Th\'eophile Vandermonde, publica el art\'iculo \textit{Remarques sur les problemes de situation} (Comentarios sobre problemas de posiciones), en el cual incluye nudos y trenzas como elementos de la geometr\'ia de posici\'on \cite{Van}.

\ 

Carl Friedrich Gauss es el siguiente matem\'atico que aporta trabajos de gran importancia para el desarrollo de la teor\'ia. De hecho, en uno de sus apuntes de 1794, titulado \textit{Una colecci\'on de nudos}, incluye 13 dibujos de nudos y m\'etodos para codificarlos. Los trabajos de Gauss fueron los que principalmente motivaron a James Clerk Maxwell quien en su libro \textit{A treatise on electricity and magnetism} (Un tratado sobre electricidad y magnetismo) incluye dos curvas cerradas que no pueden ser separadas, pero para las cuales el valor de la integral de Gauss es cero. Finalmente el a\~no 1847 represent\'o grandes avances en el desarrollo de la topolog\'ia (y por tanto de la teor\'ia de nudos) principalmente debido a que Johan Benedict Listing, un estudiante de Gauss, public\'o su monograf\'ia \textit{Vorstudien zur Topologie} (Estudios preliminares en topolog\'ia), en la cual dedica mucho espacio a los nudos. Listing fue el primero en usar la palabra topolog\'ia para referirse a la geometr\'ia de posici\'on.

\ 

En el a\~no 1876, el f\'isico escoc\'es  Sir William Thomson (Lord Kelvin) ide\'o la \textit{teor\'ia de los \'atomos de v\'ortices}. Kelvin cre\'ia que los \'atomos de la materia no eran mas que anillos de v\'ortice que se enlazaban formando distintos nudos en una sustancia llamada \'eter. La motivaci\'on de Kelvin provino del trabajo de su amigo y colega Peter G. Tait, quien hab\'ia creado un m\'etodo para producir anillos de humo. Tait por su parte, lleg\'o a estar interesado en este tipo de estudios gracias a la traducci\'on de un trabajo realizado por el matem\'atico alem\'an Hermann
von Helmholtz \textit{\"Ueber Integrale der hydrodynamischen Gleichungen, welche den Wirbelbewegungen entsprechen} (Sobre integrales de las ecuaciones hidrodinámicas que corresponden a movimientos de vórtice) en el cual aplica hidrodin\'amica a fen\'omenos electromagn\'eticos. En ese momento, crear una tabla de los elementos era de una importancia significante. La teor\'ia de Kelvin implicaba que si se creaban dos anillos de v\'ortice, enlazados el uno con el otro, entonces esto formar\'ia un \'atomo indestructible. Tait, queriendo crear un tabla de los elementos, cre\'o una clasificaci\'on de estos nudos, en lo que result\'o ser uno de los trabajos m\'as importantes en toda la historia de la teor\'ia de nudos. Como es bien sabido, la teor\'ia de Kelvin fue m\'as tarde descartada por la teor\'ia de la relatividad. De cualquier modo, es impresionante ver como a partir de una hip\'otesis refutada se gener\'o un campo formal de investigaci\'on alrededor de los nudos. 

\ 

Como se mencion\'o al inicio de esta secci\'on, la historia de la teor\'ia de nudos comprende muchos aspectos y ser\'ia imposible incluirlos todos en este art\'iculo. El lector es referido a los trabajos de J\'ozef H. Przytycki \cite{Prz, Prz2, Prz4}, que contienen gran cantidad de informaci\'on al respecto. Adem\'as, el libro \cite{PBIMW} recientemente publicado, ofrece una extensa discuci\'on hist\'orica de la teor\'ia de nudos y temas actuales de investigaci\'on incluyendo coloraciones de Fox, homolog\'ia de Yang-Baxter, polinomio de Jones, determinantes de Gram, m\'odulos de madeja, entre otros.

\section{Introducci\'on a la Teor\'ia de Nudos}\label{INTROToN}

La teor\'ia cl\'asica de nudos estudia las incrustaciones de c\'irculos m\'odulo transformaciones naturales en el espacio $\mathbb{R}^3$. Uno de los problemas fundamentales (e hist\'oricos) en la teor\'ia es la clasificaci\'on. Esta clasificaci\'on es hecha m\'odulo movimientos naturales en el espacio lo cual es llamado una isotop\'ia ambiental. 

\

Durante el siglo XIX el estudio y entendimiento de los nudos estaba lejos de ser formal. De hecho, simplemente se entend\'ia un nudo como una curva cerrada en el espacio, m\'odulo deformaciones naturales. Esta idea no estaba formalmente definida y era descrita como \textquotedblleft movimiento en el espacio sin cortar ni pegar\textquotedblright. A pesar de la falta de formalidad, con mucha intuici\'on experimental varios cient\'ificos lograron hacer aportes significativos a la teor\'ia, como por ejemplo P. G. Tait, T. P. Kirkman, C. N. Little y M. G. Haseman. En todo caso, m\'etodos precisos para diferenciar nudos no exist\'ian al momento. En esta secci\'on introducimos los conceptos de enlace poligonal, movimientos de Reidemeister, y $\Delta$-movimiento, los cuales juegan un papel fundamental en la formalizaci\'on de la teor\'ia de nudos.

\begin{definition} \
	\begin{itemize}
		\item [(a)] Un \textbf{enlace poligonal} es una colecci\'on de curvas poligonales simples disjuntas en $\mathbb{R}^{3}$. Un enlace poligonal que consiste de s\'olo una curva simple es llamado un \textbf{nudo poligonal}.
		\item [(b)] Sea $u$ un segmento o eje en un enlace poligonal $L$ en $\mathbb{R}^{3}$. Sea $\Delta$ un tri\'angulo en $\mathbb{R}^{3}$ cuya frontera consiste de tres segmentos lineales denotados por $u$, $v$ y $w$, de tal modo que $\Delta \cap L = u$. La curva poligonal $L'$, definida como $L'=(L-u)\cup  v \cup  w$, es un nuevo enlace poligonal en $\mathbb{R}^{3}$. As\'i, decimos que $L'$ ha sido obtenido desde $L$ por un $\boldsymbol{\Delta}$\textbf{-movimiento} ($\boldsymbol{\Delta}$\textbf{-move}). De forma opuesta, decimos que $L$ se obtiene a partir de $L'$ por un $\boldsymbol{\Delta^{-1}}$\textbf{-movimiento} ($\boldsymbol{\Delta^{-1}}$\textbf{-move}). Ver Figura \ref{fig:delta-move-general}. 
		\begin{figure}[ht]
			\centering
			\includegraphics[width=0.5\linewidth]{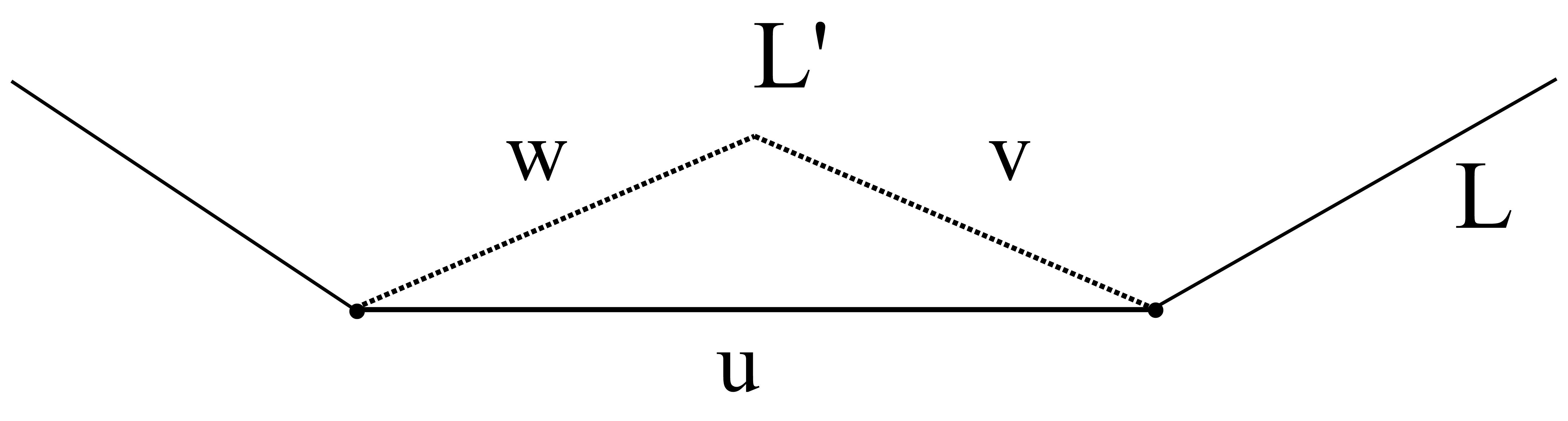}
			\caption{Se reemplaza el eje $u$ con dos ejes $v$ y $w$, a trav\'es de un $\Delta$-movimiento.} 
			\label{fig:delta-move-general}
		\end{figure}
		
		Adem\'as es permitido que el tri\'angulo $\Delta$ sea degenerado de modo que el v\'ertice $v \cap w$ est\'a en el mismo lado de $u$. En otras palabras, se permite subdivisi\'on del segmento $u$ como se ilustra en la Figura \ref{fig:delta-move-subdivision}.
		
		\begin{figure}[ht]
			\centering
			\includegraphics[width=0.5\linewidth]{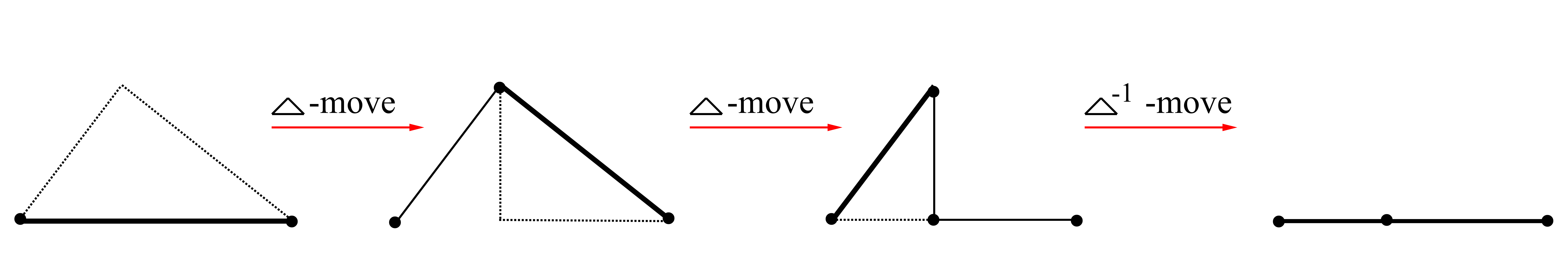}
			\caption{Subdivisi\'on como combinaci\'on de tres $\Delta$-movimientos no degenerados.}
			\label{fig:delta-move-subdivision}
		\end{figure}

		\item [(c)] Dos enlaces poligonales se dicen ser $\boldsymbol{\Delta}$\textbf{-equivalentes} (or equivalentes combinatorialmente) si uno puede ser obtenido del otro a trav\'es de una sucesi\'on finita de $\Delta$- y $\Delta^{-1}$-movimientos.
		
	\end{itemize}
	
\end{definition}

Los enlaces poligonales son usualmente presentados con diagramas planos. En palabras generales, los diagramas son definidos usando proyecciones del enlace que permiten describir la informaci\'on local en cada cruce como ``sobre'' y ``bajo''. Con miras a que esta idea est\'e bien definidia, restringimos la discusi\'on a proyecciones regulares.  Sea $p : \mathbb{R}^{3} \to \mathbb{R}^{2}$ una proyecci\'on y sea $L \subset \mathbb{R}^{3}$ un enlace. Un punto $P \in p(L)$ se dice un punto m\'ultiple de $p$ si $p^{-1}(P)$ contiene m\'as de un punto (al n\'umero de puntos en $p^{-1}(P)$ se le denomina la multiplicidad de $P$).

\begin{definition}
	La proyecci\'on $p$ se dice \textbf{regular} si
	\begin{enumerate}
		\item [(1)] $p$ tiene s\'olo un numero finito de puntos m\'ultiples y estos tienen multiplicidad dos,
		\item [(2)] ning\'un v\'ertice del enlace poligonal es preimagen de un punto m\'ultiple de $p$.
	\end{enumerate}
	
	As\'i, en el caso de una proyecci\'on regular, las siguientes situaciones ilustradas en la Figura \ref{fig:not-allowed} no son permitidas.
	\begin{figure}[ht]
		\centering
		\includegraphics[width=0.35\linewidth]{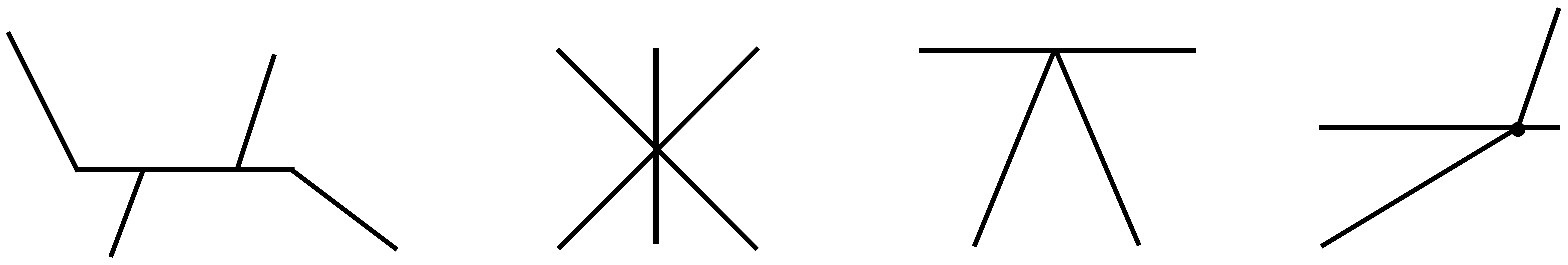}
		\caption{Situaciones no permitidas en una proyecci\'on regular de un enlace poligonal.}
		\label{fig:not-allowed}
	\end{figure}
\end{definition}	

\begin{definition}
Un \textbf{diagrama} de un enlace es una proyecci\'on regular del enlace con la informaci\'on de \textquotedblleft sobre'' y ``bajo\textquotedblright\ en cada cruce.
\end{definition}
 
Es importante observar que por razones pr\'acticas, cuando se dibuja un diagrama se asume que hay un n\'umero ``bastante grande'' de segmentos lineales que componen el nudo o enlace de modo que la traza parece estar ``curvada''.

\

En el a\~no 1927 el matem\'atico alem\'an Kurt Reidemeister demostr\'o que dos diagramas, posiblemente orientados, son isot\'opicos si y s\'olo si est\'an conectados por una sucesi\'on finita de ciertos movimientos e isotop\'ia planar. Actualmente estos movimientos son conocidos como los \textbf{movimientos de Reidemeister} aunque fueron primeramente ideados por el matem\'atico escoc\'es James Maxwell en el a\~no 1870 y demostrados independientemente por James Alexander y Garland Briggs en 1927 \cite{AB}.
 
 \begin{theorem}[Teorema de Reidemeister]\cite{Rei}
 	Dos diagramas representan enlaces $\Delta$-equivalentes si y s\'olo si los diagramas est\'an relacionados a trav\'es de una sucesi\'on finita de movimientos de Reidemeister $R_{i}$ para $i=1, 2, 3$ y deformaci\'on del plano del diagrama.
 \end{theorem}
 
\begin{figure}[h]
	\centering
	\includegraphics[width=0.7\linewidth]{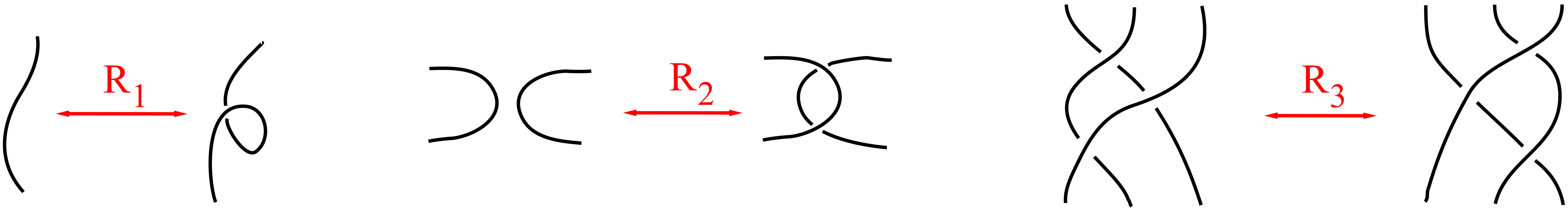}
	\caption{Movimientos de Reidemeister $R_{1}$, $R_{2}$ y $R_{3}$.}
	\label{fig:reidmoves}
\end{figure}

  La Figura \ref{fig:reidmoves} ilustra los movimientos de Reidemeister denotados por $R_{1}$, $R_{2}$ y $R_{3}$. As\'i, para distinguir nudos se usan las \textbf{invariantes}. En palabras simples, una invariante es una propiedad del enlace que se mantiene intacta bajo isotop\'ia ambiental. Para demostrar que una propiedad $P(D)$ de un diagrama $D$ es una invariante, se verifica que esta propiedad se mantenga intacta bajo los movimientos de Reidemeister. Por ejemplo, las coloraciones de Fox, el polinomio de Alexander y el polinomio de Jones son invariantes \cite{PBIMW}.

\section{Homolog\'ia de Khovanov desde el Polinomio Bracket}\label{HOM}

El matem\'atico estadounidense George David Birkhoff (1884-1944) introdujo el polinomio crom\'atico en el a\~no 1912. En el momento \'el estaba tratando de resolver el \textit{problema de los cuatro colores}. Sin entrar en detalle, el polinomio crom\'atico cuenta el n\'umero de coloraciones de los v\'ertices de un grafo, de tal manera que v\'ertices conectados por una arista son asignados diferentes colores. El primer polinomio en la teor\'ia de nudos es el polinomio de Alexander. Este importante aporte fue hecho por el matem\'atico estadounidense James Wadell Alexander (1888-1971). Se cree que Alexander sab\'ia del polinomio crom\'atico y esto lo motiv\'o a crear algo similar para nudos. Adem\'as de sus grandes aportes en ciencia, Alexander tuvo una vida muy interesante; era millonario y un gran escalador, llegando a completar ascensos populares como los Alpes Suizos y las monta\~nas rocosas en Colorado. Para m\'as detalles acerca del polinomio de Alexander, ver por ejemplo \cite{PBIMW}. 

\ 

En el a\~no 1984, un nuevo polinomio fue anunciado por el matem\'atico neozeland\'es Vaughan Frederick Randal Jones (1952-2020) \cite{Jon}. Tan importante fue este descubrimiento, que revolucion\'o la investigaci\'on en la teor\'ia de nudos a tal punto que la mayor\'ia de temas hoy en d\'ia est\'an relacionados de alguna forma al polinomio. En Agosto del a\~no 1985, Louis H. Kauffman, un matem\'atico estadounidense nacido en 1945 y uno de los m\'as activos actualmente en la investigaci\'on en teor\'ia de nudos, anunci\'o una invariante ahora conocida como el polinomio bracket de dos variables. Inicialmente Kauffman crey\'o que hab\'ia descubierto o constru\'ido una invariante completamente nueva, pero r\'apidamente se percat\'o que en realidad hab\'ia descubierto una nueva forma, m\'as sencilla, de obtener el polinomio de Jones. En ingl\'es este polinomio es conocido como el \textit{Kauffman bracket polynomial} \cite{Kau}.

\begin{definition}\label{KBPdefi}
    El polinomio bracket \textbf{cl\'asico} y \textbf{reducido}, es una funci\'on definida en el conjunto de enlaces no orientados $\mathcal{D}$ con imagen en el conjunto de los polinomios de Laurent con coeficientes enteros en la variable $A$, $\left\langle \ \right\rangle: \mathcal{D} \longrightarrow \ \mathbb{Z}[A^{\pm1}]$. Se caracteriza con las condiciones iniciales  $\left\langle \bigcirc \right\rangle=1$,  $\left\langle D \sqcup \bigcirc \right\rangle= (-A^{2}-A^{-2})\left\langle D \right\rangle$ y la relaci\'on de madeja:
		$$\left\langle\KPB\right\rangle=
		A\left\langle\KPC\right\rangle + A^{-1} \left\langle \KPD \right\rangle.$$
\end{definition}

\begin{definition}\label{Kauffmanstatedefi}
    Considere $D$ un enlace no orientado y denote por $cr(D)$ el conjunto de sus cruces. Definimos un \textbf{estado de Kauffman} (o simplemente \textbf{estado}), denotado $s$, del diagrama $D$, como una funci\'on $s: cr(D)\longrightarrow \{A, B\}$. Esta funci\'on es entendida como una asignaci\'on de un \textbf{marcador} $A$ o $B$ a \textbf{todos} los cruces como se indica en la Figura \ref{fig:kaufstatefunction}. Denotamos por $KS$ el conjunto de todos los estados. Adem\'as, en la figura se observa como cada marcador genera una \textbf{soluci\'on} (o \textbf{smoothing}) del cruce.	
    
\begin{figure}[H]
\centering
\includegraphics[width=0.9\linewidth]{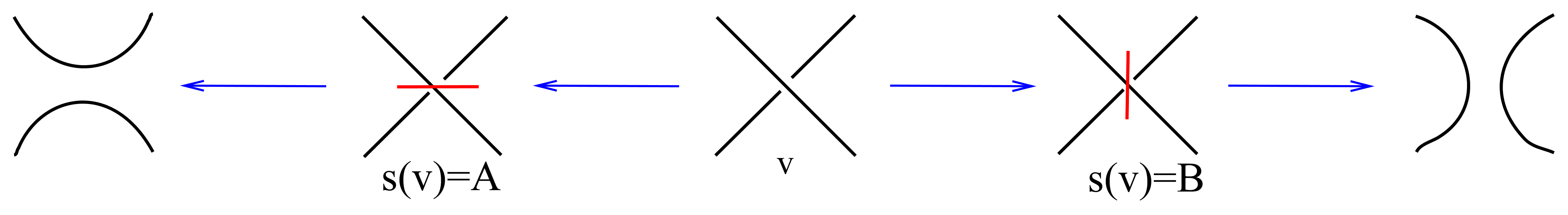}
\caption{Marcadores y smoothings de un cruce $v$.}
\label{fig:kaufstatefunction}
\end{figure}
\end{definition}

De esta manera, el polinomio bracket del diagrama $D$, est\'a dado por la siguiente f\'ormula:
	\begin{equation}
		<D>=\sum_{s \ \in \ KS}^{}A^{\mid s^{-1}(A) \mid - \mid s^{-1}(B) \mid }(-A^{2}-A^{-2})^{\mid D_{s} \mid -1}\label{EqKBP},
	\end{equation} 

donde $D_{s}$ denota el sistema de c\'irculos obtenidos por la soluci\'on de todos los cruces de acuerdo con los marcadores en el estado $s$ y $|D_{s}|$ denota el n\'umero de c\'irculos en el sistema. 

\ 

El lector pudo haber observado la palabra reducido en la Definici\'on \ref{KBPdefi}. Esto sugiere correctamente que hay varias versiones de este polinomio. Dado que m\'as adelante usaremos la versi\'on no reducida, la definimos a continuaci\'on. 

\begin{definition}
El polinomio bracket \textbf{no reducido}, es una versi\'on del polinomio bracket en la cual el enlace vac\'io $\emptyset$, es asignado el polinomio $1$. Se denota por $[ \ ]$ y con esta notaci\'on tenemos $[\emptyset]=1$, $[\bigcirc]=(-A^{2}-A^{-2})$ y $[D]=(-A^{2}-A^{-2})<D>$.
 			
\end{definition}
 
As\'i, la Ecuaci\'on \ref{EqKBP} se convierte en:

\begin{equation}\label{UnredFormula}
		[D]=\sum_{s \ \in \ KS}^{}A^{\mid s^{-1}(A) \mid - \mid s^{-1}(B) \mid }(-A^{2}-A^{-2})^{\mid D_{s} \mid}.
	\end{equation}

\begin{example}\label{KBPexampletrefoil} En este ejemplo presentamos c\'omo calcular el polinomio bracket reducido y no reducido del nudo conocido como el tr\'ebol derecho (o \textit{right-trefoil knot}). Es claro que, aunque no es la forma \'optima, el polinomio puede ser calculado directamente usando la relaci\'on de madeja en la Definici\'on \ref{KBPdefi}. Aqu\'i haremos el c\'alculo usando los estados de Kauffman. La Figura \ref{fig:examplekbp} muestra un diagrama del nudo trefoil con sus cruces enumerados \textcolor{blue}{$1,2,3$} y adem\'as contiene todos sus estados con los cruces resueltos. Estos estados est\'an nombrados por tres letras que representan el marcador asignado a cada cruce. Por ejemplo el nombre \textcolor{blue}{$ABA$} indica que los cruces \textcolor{blue}{$1$} y \textcolor{blue}{$3$} son dados un marcador \textcolor{blue}{$A$} y el cruce \textcolor{blue}{$2$} es dado un marcador \textcolor{blue}{$B$}.
		
\begin{figure}[h]
\centering
\includegraphics[width=0.55\linewidth]{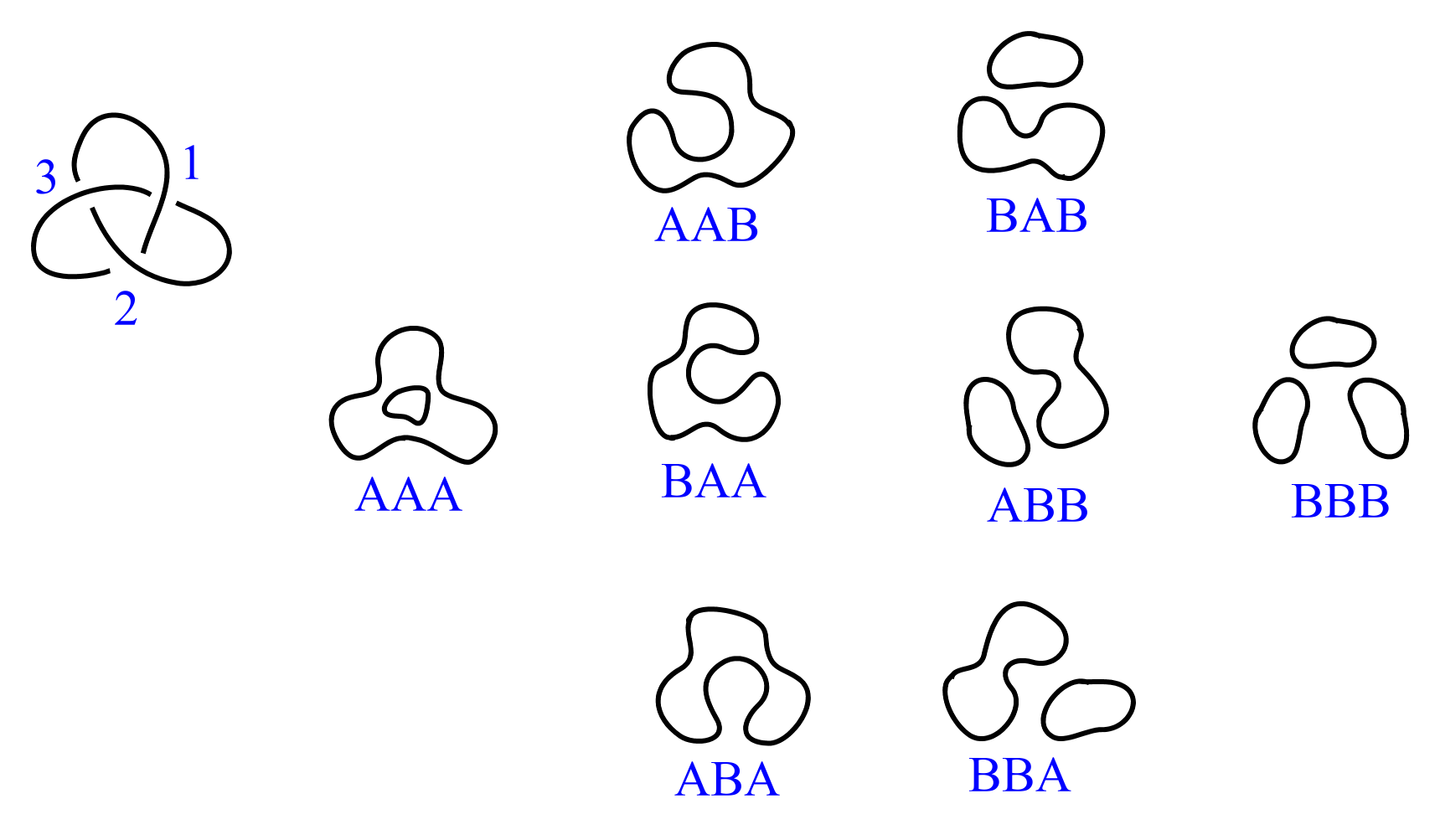}
\caption{Estados del nudo right-trefoil.}
\label{fig:examplekbp}
\end{figure}

Usando las ecuaciones previamente obtenidas, sabemos que aporta cada estado al polinomio:

	$$ \textcolor{blue}{AAA} \longrightarrow A^{3-0}(-A^2-A^{-2})^{(2-1)}=A^3(-A^2-A^{-2}),$$
	$$ \textcolor{blue}{AAB}, \textcolor{blue}{BAA}, \textcolor{blue}{ABA} \longrightarrow A^{2-1}(-A^2-A^{-2})^{1-1}=A,$$
	$$ \textcolor{blue}{BAB}, \textcolor{blue}{ABB}, \textcolor{blue}{BBA} \longrightarrow A^{1-2}(A^{2}-A^{-2})^{2-1}=A^{-1}(-A^{2}-A^{-2}), \ y   $$
	$$\textcolor{blue}{BBB} \longrightarrow A^{0-3}(-A^{2}-A^{-2})^{3-1}=A^{-3}(-A^{2}-A^{-2})^{2}.$$
	
	De esta forma, el polinomio bracket reducido del nudo trefoil esta dado por:
 $$A^3(-A^{2}-A^{-2})+3A+3A^{-1}(-A^{2}-A^{-2})+A^{-3}(-A^{2}-A^{-2})^{2} $$ 
	$$=A^{-7}-A^{5}-A^{-3}.$$
 
	Adicionalmente, la versi\'on no reducida est\'a dada por:
 
	$$(-A^{2}-A^{-2})(A^{-7}-A^{5}-A^{-3})=-A^{-9}+A^{7}+A^{3}+A^{-1}.$$
	
\end{example}

Observe que los terminos en la Ecuaci\'on \ref{UnredFormula} provienen de los estados del diagrama y por tanto tienen una interpretaci\'on geom\'etrica. Adem\'as, esta f\'ormula permite tener una correspondencia 1-1 entre los c\'irculos en $D_{s}$ y los factores $(-A^{2}-A^{-2})$. Esto motiva la siguiente definici\'on.

\begin{definition}\label{enhancedstatesdefi}
	Un \textbf{estado de Kauffman mejorado} $S$ del diagrama $D$, es un estado de Kauffman $s$  junto a una funci\'on  $\varepsilon: D_{s}  \longrightarrow \{+, -\}$, la cual asigna a cada c\'irculo de $D_{s}$ un signo positivo o un signo negativo.
\end{definition}	

Se denota al conjunto de los estados mejorados de Kauffman por $EKS$ (por sus siglas en ingl\'es \textit{Enhanced Kauffman States}). Puesto que cada c\'irculo en $D_{s}$ es dado un signo positivo o negativo, se deduce que para cada estado $s$ existen $2^{\mid D_{s} \mid}$ estados mejorados. As\'i, el polinomio bracket puede ser expresado como una suma de monomios, provenientes de los estados mejorados, de la siguiente manera: 

\begin{equation}\label{monosum}
	[D]=\sum_{S \ \in \ EKS}^{}(-1)^{\mid D_{s} \mid}A^{\mid s^{-1}(A) \mid - \mid s^{-1}(B) \mid }(A^{2})^{\mid \varepsilon^{-1}(+) \mid - \mid \varepsilon^{-1}(-) \mid }.
\end{equation}

Denote por $\sigma (s)$ la diferencia entre el n\'umero de marcadores $A$ (tambi\'en llamados marcadores positivos) y marcadores $B$ (marcadores negativos) en el estado. Denote por $\tau(S)$ la diferencia entre el n\'umero de signos positivos y negativos en el estado mejorado. Esto es: 
$$\sigma (s)=|s^{-1}(A)|-|s^{-1}(B)| \ \ \ \ \ \ \ y \ \ \ \ \ \ \ \tau (S)=|\varepsilon^{-1}(+)|-|\varepsilon^{-1}(-)|.$$

De esta forma, la Ecuaci\'on \ref{monosum} se modifica a:
\begin{equation}\label{statesumunreduced}
	[D]=\sum_{S \ \in \ EKS}^{}(-1)^{\mid D_{s} \mid}A^{\sigma (s)+2\tau (S)},
\end{equation}			
la cual se conoce como la f\'ormula para el polinomio bracket no reducido en t\'erminos de los estados mejorados. Estos estados mejorados forman una base para los grupos de cadena y el complejo de cadena de Khovanov definidos a continuaci\'on.

\begin{definition}\ \label{basis}
	
	\begin{enumerate}
		\item [(i)] El \textbf{bigrado} en los estados mejorados se define con el siguiente conjunto: $$\mathcal{S}_{a,b}(D)=\mathcal{S}_{a,b}=\left\lbrace  S \ \in \ EKS \ \mid \  a=\sigma (s), \  b=\sigma (s)+2\tau (S) \right\rbrace. $$
		\item [(ii)] Los \textbf{grupos de cadena} $\mathcal{C}_{a,b}(D)=\mathcal{C}_{a,b}$, se definen como los grupos abelianos libres con base $\mathcal{S}_{a,b}(D)=\mathcal{S}_{a,b}$, i.e. $\mathcal{C}_{a,b}=\mathbb{Z}\mathcal{S}_{a,b}$. As\'i, $\mathcal{C}(D)=\displaystyle\bigoplus_{a,b \ \in \ \mathbb{Z}}\mathcal{C}_{a,b}(D)$ es grupo libre abeliano bigradado.
		\item [(iii)] Para un diagrama $D$ de un enlace, se define el \textbf{complejo de cadena} $\mathscr{C}(D)= \left\lbrace \left( C_{a,b}, \partial_{a,b} \right) \right\rbrace  $,
		donde 
		el \textbf{mapa diferencial} $\partial_{a,b}: \mathcal{C}_{a,b} \longrightarrow \mathcal{C}_{a-2,b}$ es definido por:
		$$\partial_{a,b}(S)= \sum_{S' \ \in \ \mathcal{S}_{a-2,b} }^{}(-1)^{t(S,S')}(S,S')S'.$$
	\end{enumerate}
	
\end{definition}

En el item (iii) de la definici\'on previa, $S \in \mathcal{S}_{a,b}$ y $(S,S')$ es lo que se conoce como el n\'umero incidental entre $S$ y $S'$. Tenemos que $(S,S')=0$ o $(S, S')=1$. Adem\'as, sabemos que es igual a $1$ si y s\'olo si las siguientes condiciones se cumplen:
 \begin{enumerate}
 	\item $S$ y $S'$ son id\'enticos, excepto en un cruce, denotado por $v$. Adem\'as, $s(v)=A$ y $s'(v)=B$, donde $S$ y $S'$ son estados mejorados de $s$ y $s'$ respectivamente. 
 	\item $\tau (S')=\tau (S)+1$ y cada componente de $D_{s}$ que no interact\'ua con el cruce $v$ mantiene su signo por $D_{s'}$. 
 \end{enumerate}	
Note que la condicion (1) describe el hecho que el valor de $\sigma$ decrece en $2$, mientras que la condici\'on (2) representa el hecho de que el n\'umero de signos negativos decrece o el n\'umero de cruces positivos incrementa. Finalmente, el n\'umero $(-1)^{t(S,S')}$ requiere un orden en los cruces del diagrama $D$. El n\'umero $t(S,S')$ es definido como el n\'umero de cruces en $D$ asignados un marcador $B$ en $S$ despu\'es del cruce $v$ en el orden estipulado. Este orden es suficiente para que el mapa diferencial cumpla con la condici\'on $\partial_{a-2,b}\circ \partial_{a,b}=0$. De cualquier modo, es importante mencionar que la homolog\'ia no depende del orden de los cruces \cite{Kho}.

\ 

\begin{figure}[H]
	\centering
	\includegraphics[width=0.85\linewidth]{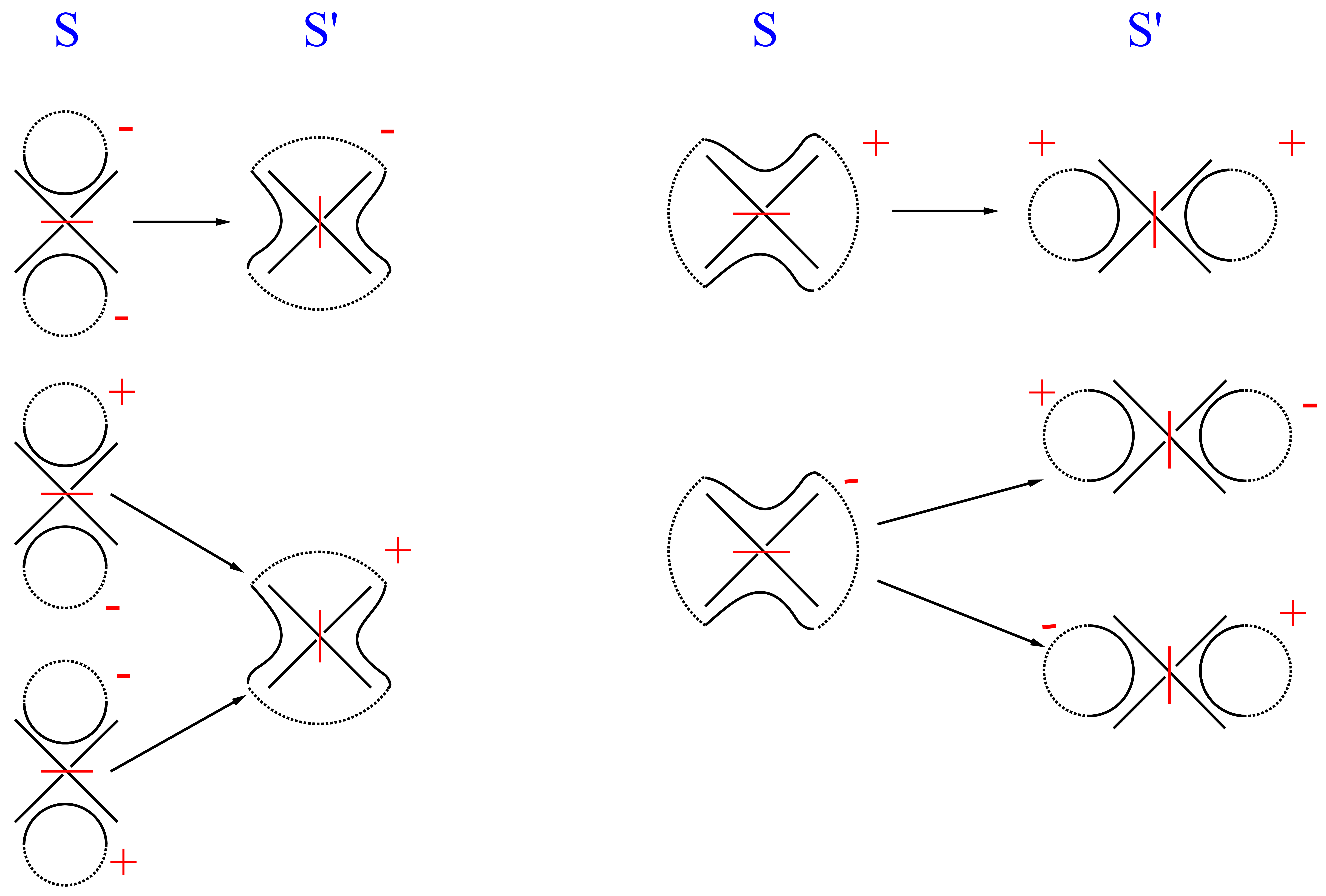}
	\caption{Signos de los c\'irculos despu\'es de uni\'on (izquierda) y separaci\'on (derecha) de los componentes.}
	\label{fig:cases-merged}
\end{figure}

Con la notaci\'on dada hasta el momento, podemos presentar la Figura \ref{fig:cases-merged} que muestra los posibles escenarios en que los estados mejorados $S$ y $S'$ son incidentales ( $(S,S')=1$). Por ejemplo, cuando dos c\'irculos con signos diferentes son ``unidos'', el resultado debe ser un c\'irculo con un signo positivo.

	\begin{definition}\label{khdefi}
	La \textbf{homolog\'ia de Khovanov} del diagrama $D$ se define como la homolog\'ia del complejo de cadena $\mathscr{C}(D)$:		
	$$H_{a,b}(D)=\dfrac{ker (\partial_{a,b}) }{im(\partial_{a+2,b})}.$$
\end{definition}

Por supuesto que la propiedad m\'as importante de la homolog\'ia de Khovanov es su invarianza bajo los movimientos de Reidemeister de segundo y tercer tipo. El siguiente teorema establece el resultado; ver \cite{Vir1, Vir2} para una demostraci\'on.

\begin{theorem}\label{khinvariance}
	Sea $D$ un diagrama de un enlace. Los grupos de homolog\'ia
	$$H_{a,b}(D)=\dfrac{ker (\partial_{a,b}) }{im(\partial_{a+2,b})},$$
	son invariantes bajo los movimientos de Reidemeister de segundo y tercer tipo. Por lo tanto, son invariantes de enlaces no orientados enmarcados\footnote{En ingl\'es \textit{unoriented framed knots}.}. Adem\'as, el efecto del movimiento de Reidemeister de primer tipo $R_{1}$ (positivo o negativo) est\'a dado por $H_{a,b}(R_{1+}(D))=H_{a+1,b+3}(D)$ y $H_{a,b}(R_{1-}(D))=H_{a-1,b-3}(D)$. Estos grupos categorifican el polinomio bracket no reducido y son llamados los grupos homol\'ogicos enmarcados de Khovanov.   
\end{theorem}

\subsection{Ejemplo: Homolog\'ia del nudo tr\'ebol}

\begin{figure}[h]
	\centering
	\includegraphics[width=0.7\linewidth]{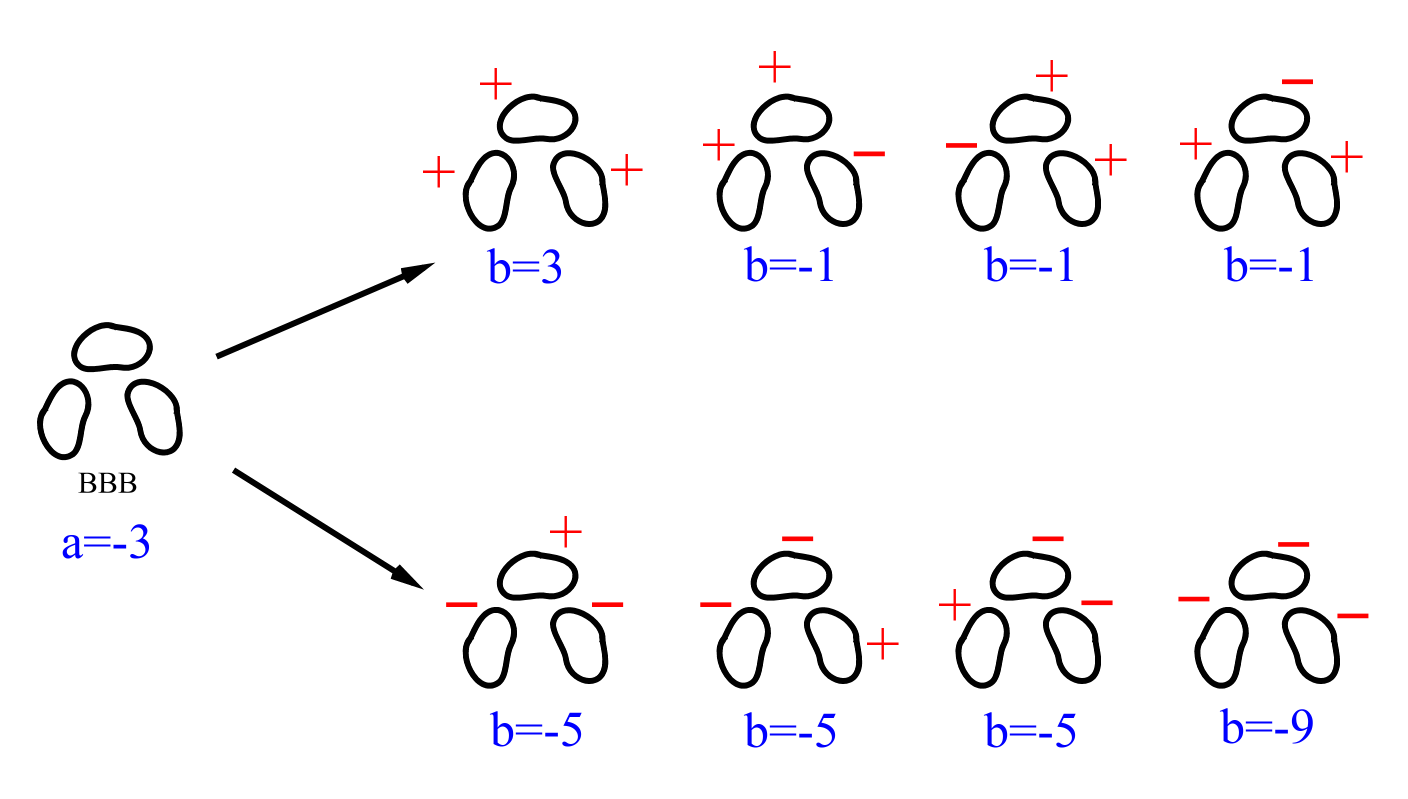}
	\caption{Estados mejorados del estado $BBB$.}
	\label{fig:enhancedbbb}
\end{figure}

Ahora presentamos a manera de ejemplo, el c\'alculo de la homolog\'ia de Khovanov para el nudo tr\'ebol. Primero se debe asignar un orden a los cruces y obtener todos los estados, tal cual se observa en la Figura \ref{fig:examplekbp}. Una vez se tienen los estados de Kauffman, es posible determinar los estados mejorados y los valores de $a=\sigma(s)$ y $b=\sigma(s)+2\tau(S)$. Por ejemplo, la Figura \ref{fig:enhancedbbb} muestra los estados mejorados que corresponden al estado $BBB$. 

\begin{figure}[h]
	\centering
	\includegraphics[width=0.7\linewidth]{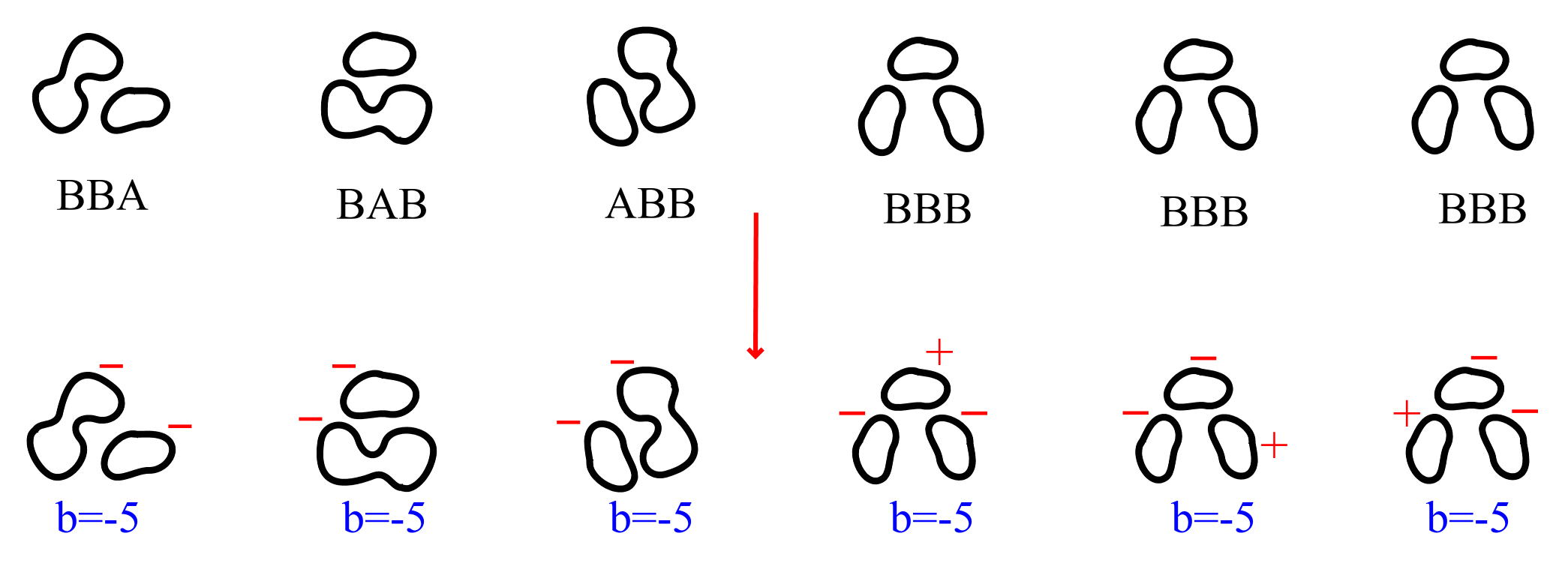}
	\caption{Estados mejorados del nudo trefoil para los cuales $b=-5$.}
	\label{fig:complexb-5}
\end{figure}

Teniendo los valores del bigrado $(a,b)$ para los estados mejorados del nudo trefoil, el pr\'oximo paso es construir el complejo de cadena. Recordamos que $\partial_{a,b}: \mathcal{C}_{a,b} \longrightarrow \mathcal{C}_{a-2,b}$ es dado en la Definici\'on \ref{basis}. Observe que el valor de $b$ no cambia, lo cual sugiere que para cada valor de $b$ un complejo es construido. La Figura \ref{fig:complexb-5} ilustra los estados mejorados para los cuales obtenemos el valor $b=\sigma(s)+2\tau(S)=-5$.

\begin{figure}[H]
	\centering
	\begin{overpic}[scale=1.4]{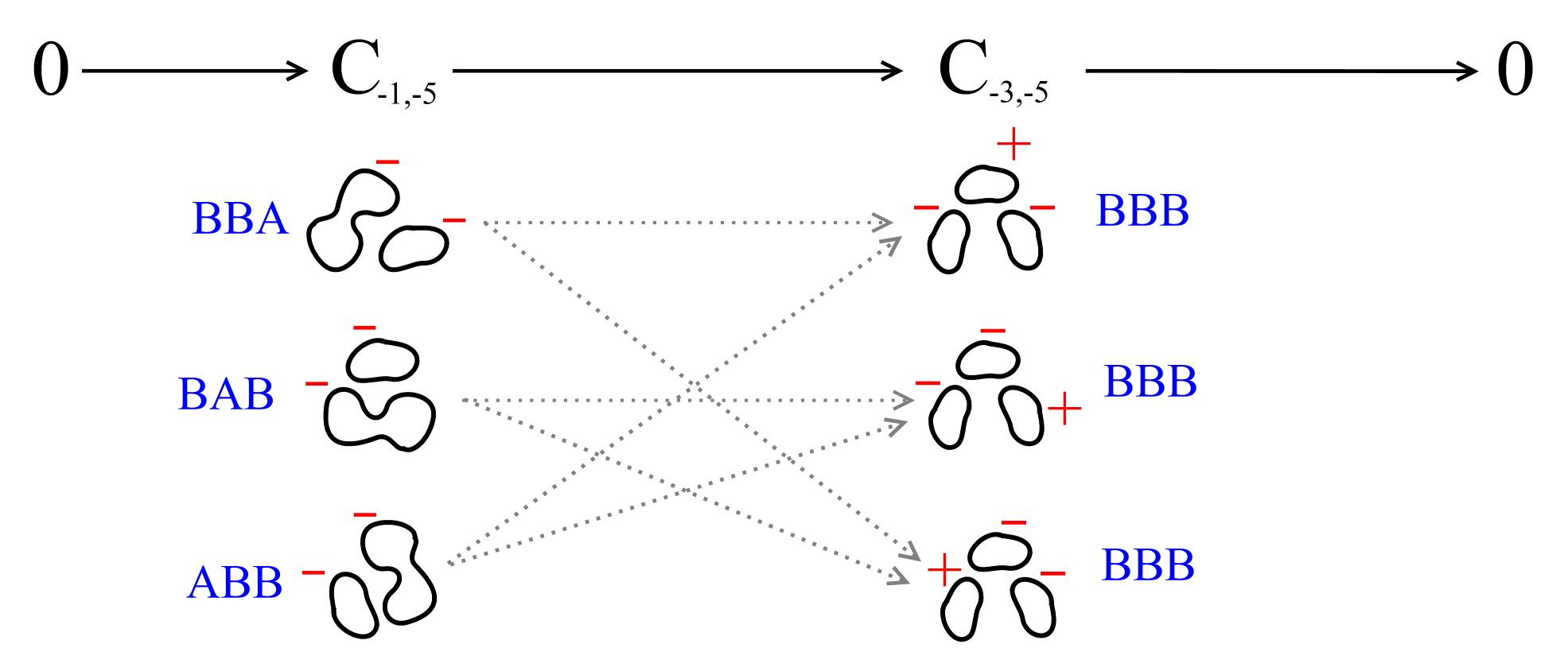}
		\put(145,160){\Large \color{blue} $ \partial_{-1,-5}$}
		\put(295,160){\Large \color{blue} $\partial_{-3,-5}$}
	\end{overpic}
	\caption{Complejo de cadena para $b=-5$.}
	\label{fig:chaincomplexb-5}
\end{figure}

De este modo tenemos que:
$\mathcal{C}_{-1,-5}=\mathbb{Z}\mathcal{S}_{-1,-5}$ donde $$\mathcal{S}_{-1,-5}=\left\lbrace  S \ \in \ EKS \ \mid \  \sigma (s)=-1, \  \sigma (s)+2\tau (S)=-5 \right\rbrace.$$
Similarmente, $\mathcal{C}_{-3,-5}=\mathbb{Z}\mathcal{S}_{-3,-5}$ donde $$\mathcal{S}_{-3,-5}=\left\lbrace  S \ \in \ EKS \ \mid \  \sigma (s)=-3, \  \sigma (s)+2\tau (S)=-5 \right\rbrace.$$ 
As\'i, el complejo obtenido para $b=-5$ es ilustrado en la Figura \ref{fig:chaincomplexb-5}. Adem\'as, la figura muestra los pares de estados incidentales.	

Se puede observar que $H_{-1,-5}=\dfrac{ker (\partial_{-1,-5}) }{im(\partial_{1,-5})}$ es trivial.\\
Ahora, para obtener $H_{-3,-5}=\dfrac{ker (\partial_{-3,-5}) }{im(\partial_{-1,-5})}$, note que $ker (\partial_{-3,-5})=\mathcal{C}_{-3,-5}$ es el grupo abeliano generado por los siguientes estados mejorados, los cuales son denotados por $x, y, z$  para simplificar el procedimiento, como se ilustra en la Figura \ref{fig:kernel-3-5}:
\begin{figure}[H]
	\centering
	\includegraphics[width=0.45\linewidth]{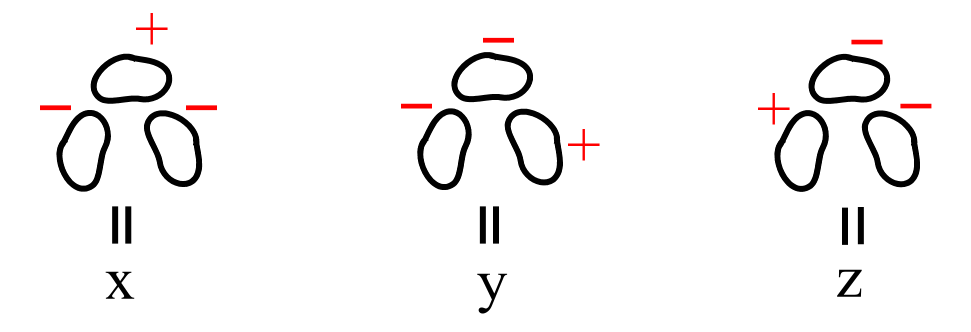}
	\caption{Generadores de $ker (\partial_{-3,-5})$.}
	\label{fig:kernel-3-5}
\end{figure}

Por otro lado, note que la imagen de $\partial_{-1,-5}$, es el grupo abeliano cuya base esta dada por:

\begin{figure}[ht]
	\centering
	\includegraphics[width=0.7\linewidth]{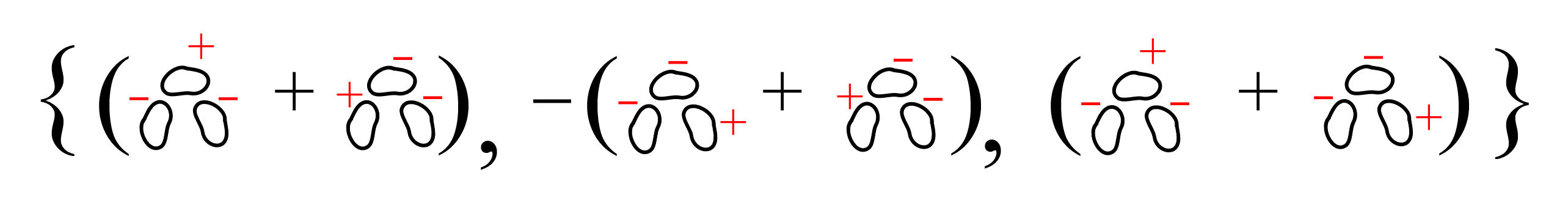}
	\caption{Base para $im (\partial_{-1,-5})$.}
	\label{fig:im-1-5}
\end{figure}

Por lo tanto, el grupo homol\'ogico $H_{-3,-5}$ esta dado por $$H_{-3,-5}=\dfrac{\mathbb{Z}\left\lbrace x, y, z\right\rbrace }{\left\lbrace x+z, -y-z, x+y\right\rbrace }.$$
Estas relaciones son equivalentes a tener $z=-x=y$ y $2x=0$. As\'i, este cociente es $\mathbb{Z}_{2}$, un grupo de torsi\'on. El proceso para los otros grupos homol\'ogicos es an\'alogo (ver por ejemplo \cite{MV}). La Tabla \ref{KH of the right-trefoil knot in Viro} contiene la homolog\'ia de Khovanov completa para el nudo trefoil.

\begin{table}[H]
	\centering
	\begin{tabular}{c|||c|c|c|c|c} 
		\ \ $\textbf{\textcolor{red}{b}} \ \ | \  \ \textbf{\textcolor{red}{a}}$ \ \ & \ \ \ -3 \ \ \       & \ \ \ -1 \ \ \ & \ \ \ 1 \ \ \ & \ \ \ 3 \ \ \ \\
		\hline \hline 
		7                                           &                      &             &           & $\mathbb{Z}$  \\ 
		\hline
		3                                           &                     &             &            & $\mathbb{Z}$    \\
		\hline
		-1                                          &                     &$\mathbb{Z}$ &            &    \\
		\hline
		-5                                          &    $\mathbb{Z}_{2}$ &             &             & \\
		\hline
		-9                                          &  $\mathbb{Z}$       &              &             &    \\
		\hline
		
	\end{tabular} 
	\caption{Homolog\'ia de Khovanov del nudo trefoil.}\label{KH of the right-trefoil knot in Viro}
\end{table}

\section{Sucesi\'on Larga de Khovanov}\label{suclarga}

De una forma similar en la que se construye la homolog\'ia de Khovanov, podemos categorificar la relaci\'on de madeja en la definici\'on del polinomio bracket (Definici\'on \ref{KBPdefi}). Recuerde que en la Definici\'on \ref{basis} se da el siguiente conjunto:
$$\mathcal{S}_{a,b}(D)=\mathcal{S}_{a,b}=\left\lbrace  S \ \in \ EKS \ \mid \  a=\sigma (s), \  b=\sigma (s)+2\tau (S) \right\rbrace,$$
el cual sirve como base para los grupos abelianos libres  $\mathcal{C}_{a,b}=\mathbb{Z}\mathcal{S}_{a,b}$. Sea $v$ un cruce fijo del diagrama $D$. Considere los conjuntos 
$\mathcal{S}_{a,b}^{A,v}$ y $\mathcal{S}_{a,b}^{B,v}$ definidos por:
	$$\mathcal{S}_{a,b}^{A,v}=\left\lbrace S \ \in \ S_{a,b} \ \mid \ s(v)=A \right\rbrace \ \ \ \ \ y  \ \ \ \  \ \mathcal{S}_{a,b}^{B,v}=\left\lbrace S \ \in \ S_{a,b} \ \mid \ s(v)=B \right\rbrace,$$
i.e. $\mathcal{S}_{a,b}^{A,v}$ contiene todos los estados mejorados $S$ con bigrado $(a,b)$ que tienen un marcador $A$ en el cruce $v$; de forma an\'aloga, $\mathcal{S}_{a,b}^{B,v}$ contiene todos los estados mejorados $S$ con bigrado $(a,b)$ que tienen un marcador $B$ en el cruce $v$. De esta forma, podemos observar que: 	
$$\mathcal{S}_{a,b} \ \ = \ \   \mathcal{S}_{a,b}^{A,v} \ \ \sqcup \ \ \mathcal{S}_{a,b}^{B,v}.$$

Usando la misma notaci\'on de la construcci\'on de la homolog\'ia de Khovanov, denotamos los grupos abelianos generados por estos conjuntos por:
$$\mathbb{Z}\mathcal{S}_{a,b}^{A,v} \ \ = \ \ \mathcal{C}_{a,b}^{\includegraphics[width=0.04\textheight]{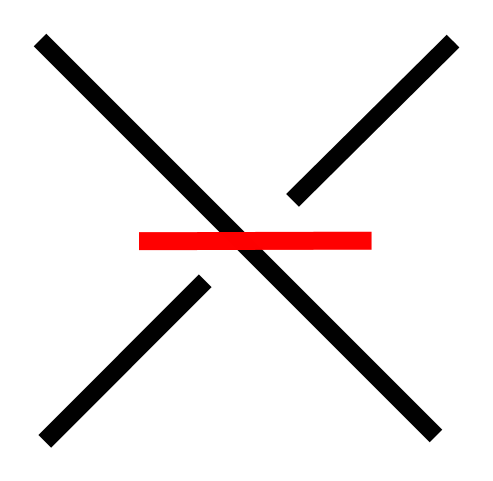}} \ \ \ \ \  y \ \ \ \ \ \mathbb{Z}\mathcal{S}_{a,b}^{B,v} \ \ = \ \ \mathcal{C}_{a,b}^{\includegraphics[width=0.04\textheight]{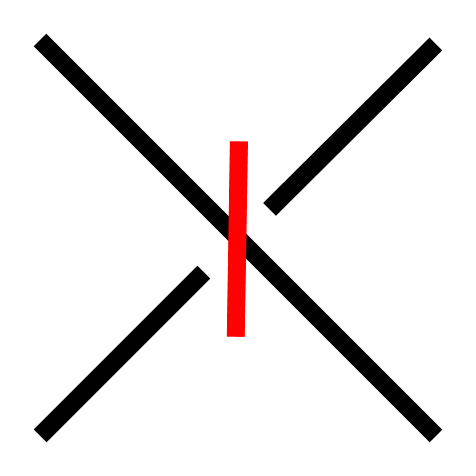}},$$
donde $v= \textbf{\KPB}$.  As\'i, al nivel de grupos se obtienes las siguientes igualdades:
$$\mathcal{C}_{a,b} \ \ = \ \ \mathcal{C}_{a,b}^{\includegraphics[width=0.04\textheight]{imagenes/positive_marker}} \ \ \oplus \ \ \mathcal{C}_{a,b}^{\includegraphics[width=0.04\textheight]{imagenes/negative_marker}} \ \ \ o \ \  equivalentemente \ \ \ \mathbb{Z}\mathcal{S}_{a,b} \ \ = \ \ \mathbb{Z}\mathcal{S}_{a,b}^{A,v} \ \ \sqcup \ \ \mathbb{Z}\mathcal{S}_{a,b}^{B,v}.$$ 

Observe que el complejo $( \mathcal{C}_{a,b}^{\includegraphics[width=0.04\textheight]{imagenes/negative_marker}}, \partial_{a,b}) $ es un sub-complejo de cadena de $\left( \mathcal{C}_{a,b}, \partial_{a,b}  \right) $, en otras palabras $\partial(\mathcal{C}_{a,b}^{\includegraphics[width=0.04\textheight]{imagenes/negative_marker}}) \subset \mathcal{C}_{a-2,b}^{\includegraphics[width=0.04\textheight]{imagenes/negative_marker}}$. En contraste, note que esto no es necesariamente cierto para $(\mathcal{C}_{a,b}^{\includegraphics[width=0.04\textheight]{imagenes/positive_marker}}, \partial_{a,b})$ puesto que $\partial_{a,b}$ podr\'ia cambiar el marcador del cruce de $A$ a $B$, como fue discutido mientras se construy\'o la homolog\'ia de Khovanov. La siguiente sucesi\'on exacta corta puede ser escrita:

$$0 \xrightarrow{\makebox[1.6cm]{}}\mathcal{C}_{a,b}^{\includegraphics[width=0.04\textheight]{imagenes/negative_marker}} \xrightarrow{\makebox[1.6cm]{\textcolor{blue}{$\alpha$}}} \mathcal{C}_{a,b}^{\includegraphics[width=0.04\textheight]{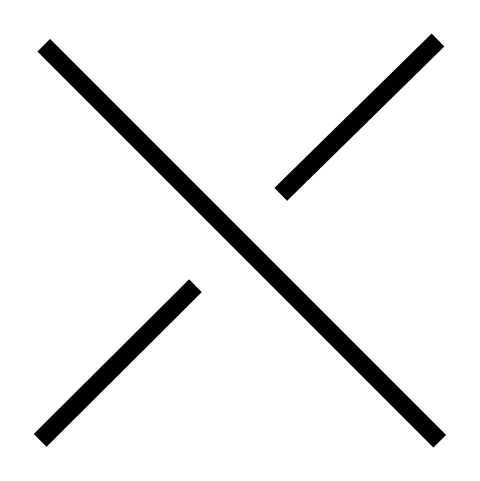}} \xrightarrow{\makebox[1.6cm]{\textcolor{blue}{$\beta$}}}  \quotient{\textstyle \mathcal{C}_{a,b}^{\includegraphics[width=0.04\textheight]{imagenes/Crossing}}}{\mathcal{C}_{a,b}^{\includegraphics[width=0.04\textheight]{imagenes/negative_marker}}}  \xrightarrow{\makebox[1.6cm]{}} 0.$$

En esta sucesi\'on, $\alpha$ toma un estado mejorado del diagrama $D$ con marcador $B$ en el cruce $v$ y lo env\'ia al estado mejorado de $D$ que asigna un marcador $B$ al cruce $v$ y los dem\'as cruces mantienen sus marcadores y los c\'irculos mantienen sus signos. Por otro lado, el mapa $\beta$ env\'ia a los estados mejorados con marcador $B$ en el cruce $v$ a cero y a cada estado mejorado con marcador $A$ en el cruce $v$ al estado mejorado de $D$ con el cruce $v$ dado un marcador $A$ y los dem\'as cruces y c\'irculos preservan sus marcadores y signos, respectivamente.\\

Observe que $\quotient{\textstyle \mathcal{C}_{a,b}^{\includegraphics[width=0.04\textheight]{imagenes/Crossing}}}{\mathcal{C}_{a,b}^{\includegraphics[width=0.04\textheight]{imagenes/negative_marker}}} \ \ = \ \ \mathcal{C}_{a,b}^{\includegraphics[width=0.04\textheight]{imagenes/positive_marker}}$ y existe un mapa conector entre los complejos de cadena de la homolog\'ia de Khovanov: 

$$\quotient{\textstyle \mathcal{C}_{a,b}^{\includegraphics[width=0.04\textheight]{imagenes/Crossing}}}{\mathcal{C}_{a,b}^{\includegraphics[width=0.04\textheight]{imagenes/negative_marker}}} \xrightarrow{\makebox[3cm]{\textcolor{blue}{$\partial^{Conn}_{a,b}$}}} \mathcal{C}_{a-2,b}^{\includegraphics[width=0.04\textheight]{imagenes/negative_marker}}.$$
As\'i, obtenemos la siguiente sucesi\'on larga de homolog\'ia:
\begin{equation*}
	\begin{split}
		\cdots \xrightarrow{\makebox[1cm]{}} H_{a,b}( \includegraphics[width=0.04\textheight]{imagenes/negative_marker}) &   \xrightarrow{\makebox[1.5cm]{\textcolor{blue}{$\alpha_{\ast}$}}} H_{a,b}(\includegraphics[width=0.04\textheight]{imagenes/Crossing})\xrightarrow{\makebox[1.8cm]{\textcolor{blue}{$\beta_{\ast}$}}}H_{a,b}\left(\quotient{\textstyle \includegraphics[width=0.04\textheight]{imagenes/Crossing}} {\includegraphics[width=0.04\textheight]{imagenes/negative_marker}}\right)  \\
		\xrightarrow{\makebox[2cm]{$  \textcolor{blue}{\left( \partial^{Conn}_{a,b} \right)_{\ast}}$}} H_{a-2,b}(\includegraphics[width=0.04\textheight]{imagenes/negative_marker}) & \xrightarrow{\makebox[1.5cm]{\textcolor{blue}{$\alpha_{\ast}$}}}H_{a-2,b}(\includegraphics[width=0.04\textheight]{imagenes/Crossing}) \xrightarrow{\makebox[1.5cm]{\textcolor{blue}{$\beta_{\ast}$}}}H_{a-2,b}\left(\quotient{\textstyle \includegraphics[width=0.04\textheight]{imagenes/Crossing}} {\includegraphics[width=0.04\textheight]{imagenes/negative_marker}}\right)  \\
		\xrightarrow{\makebox[2cm]{$ \textcolor{blue}{\left( \partial^{Conn}_{a,b} \right)_{\ast}}$}} H_{a-4,b}(\includegraphics[width=0.04\textheight]{imagenes/negative_marker}) &  \xrightarrow{\makebox[1.5cm]{\textcolor{blue}{$\alpha_{\ast}$}}}H_{a-4,b}(\includegraphics[width=0.04\textheight]{imagenes/Crossing})\xrightarrow{\makebox[1.8cm]{\textcolor{blue}{$\beta_{\ast}$}}}\cdots.  \\ 
	\end{split}
\end{equation*}	

Note las siguientes igualdades:
$$\mathcal{C}_{a,b}^{\includegraphics[width=0.04\textheight]{imagenes/negative_marker}} \ \ = \ \ \mathcal{C}_{a+1,b+1}^{\includegraphics[width=0.04\textheight]{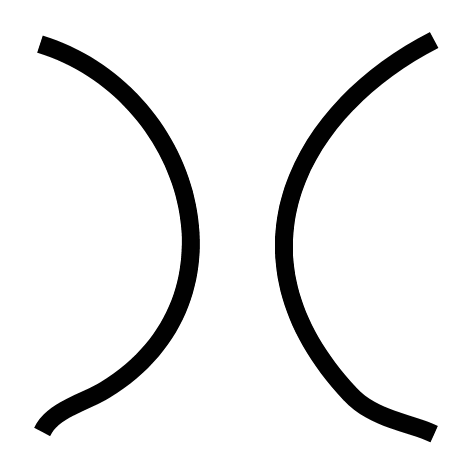}} \ \ \ \ \ \ y \ \ \ \ \ \ \quotient{\textstyle \mathcal{C}_{a,b}^{\includegraphics[width=0.04\textheight]{imagenes/Crossing}}}{\mathcal{C}_{a,b}^{\includegraphics[width=0.04\textheight]{imagenes/negative_marker}}} \ \ = \ \ \mathcal{C}_{a-1,b-1}^{\includegraphics[width=0.04\textheight]{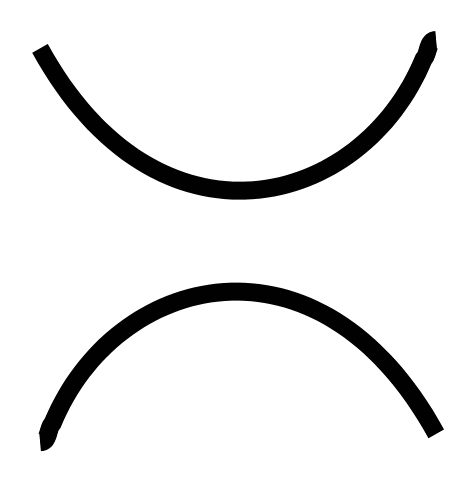}}. $$


De tal modo que, podemos obtener la siguiente sucesi\'on corta de complejo de cadenas de diagramas:
$$0 \xrightarrow{\makebox[1.5cm]{}}\mathcal{C}_{a+1,b+1}^{\includegraphics[width=0.04\textheight]{imagenes/b_smoothing}} \xrightarrow{\makebox[1.5cm]{\textcolor{blue}{$\alpha$}}} \mathcal{C}_{a,b}^{\includegraphics[width=0.04\textheight]{imagenes/Crossing}} \xrightarrow{\makebox[1.5cm]{\textcolor{blue}{$\beta$}}}\mathcal{C}_{a-1,b-1}^{\includegraphics[width=0.04\textheight]{imagenes/a_smoothing}}    \xrightarrow{\makebox[1.5cm]{}} 0,$$

la cual, genera la siguiente sucesi\'on conocida como la \textbf{sucesi\'on larga de la homolog\'ia de Khovanov}:

\begin{equation}\label{LES}
	\begin{split}
		\cdots \xrightarrow{\makebox[1cm]{}} H_{a+1,b+1}(\includegraphics[width=0.035\textheight]{imagenes/b_smoothing}) & \xrightarrow{\makebox[1cm]{\textcolor{blue}{$\alpha_{\ast}$}}} H_{a,b}(\includegraphics[width=0.035\textheight]{imagenes/Crossing})\xrightarrow{\makebox[1.3cm]{\textcolor{blue}{$\beta_{\ast}$}}}H_{a-1,b-1}(\includegraphics[width=0.035\textheight]{imagenes/a_smoothing}) \\
		\xrightarrow{\makebox[2cm]{$  \textcolor{blue}{\left( \partial^{Conn}_{a,b} \right)_{\ast}}$}} H_{a-1,b+1}(\includegraphics[width=0.035\textheight]{imagenes/b_smoothing}) & \xrightarrow{\makebox[1cm]{\textcolor{blue}{$\alpha_{\ast}$}}} H_{a-2,b}(\includegraphics[width=0.035\textheight]{imagenes/Crossing})\xrightarrow{\makebox[1cm]{\textcolor{blue}{$\beta_{\ast}$}}}H_{a-3,b-1}(\includegraphics[width=0.035\textheight]{imagenes/a_smoothing}) \\
		\xrightarrow{\makebox[2cm]{$  \textcolor{blue}{\left( \partial^{Conn}_{a,b} \right)_{\ast}}$}} H_{a-3,b+1}(\includegraphics[width=0.035\textheight]{imagenes/b_smoothing}) & \xrightarrow{\makebox[1cm]{\textcolor{blue}{$\alpha_{\ast}$}}}H_{a-4,b}(\includegraphics[width=0.035\textheight]{imagenes/Crossing})\xrightarrow{\makebox[1.3cm]{\textcolor{blue}{$\beta_{\ast}$}}} \cdots. \\
	\end{split}
\end{equation}

El siguiente resultado se obtiene directamente de la construcci\'on previa:

\begin{corollary}\    \label{Coro1}
	\begin{enumerate}
		\item [(1)] Si $H_{a+1,b+1}(\includegraphics[width=0.025\textheight]{imagenes/b_smoothing})=0$ entonces $\beta_*: H_{a+1,b+1}(\includegraphics[width=0.025\textheight]{imagenes/Crossing}) \to H_{a,b}(\includegraphics[width=0.025\textheight]{imagenes/a_smoothing})$ es un monomorfismo.
		\item[(2)] Si $H_{a-1,b+1}(\includegraphics[width=0.025\textheight]{imagenes/b_smoothing})=0$ entonces $\beta_*: H_{a+1,b+1}(\includegraphics[width=0.025\textheight]{imagenes/Crossing}) \to H_{a,b}(\includegraphics[width=0.025\textheight]{imagenes/a_smoothing})$ es un epimorfismo.
		\item [(3)] Si $H_{a+1,b+1}(\includegraphics[width=0.025\textheight]{imagenes/b_smoothing})=0=H_{a-1,b+1}(\includegraphics[width=0.025\textheight]{imagenes/b_smoothing})$ entonces $\beta_*: H_{a+1,b+1}(\includegraphics[width=0.025\textheight]{imagenes/Crossing}) \to H_{a,b}(\includegraphics[width=0.025\textheight]{imagenes/a_smoothing})$ es un isomorfismo.
	\end{enumerate}
\end{corollary}

\section{Homolog\'ia de Nudos Toroidales T(2,n)}\label{HOMTORUS}
La familia de nudos toroidales $T(2,n)$ ha siempre sido de gran inter\'es para la investigaci\'on en teor\'ia de nudos. La homolog\'ia de Khovanov para estos nudos fue calculada en primera instancia por el creador de esta teor\'ia de homolog\'ia, el matem\'atico ruso Mikhail Khovanov\footnote{Una rese\~na hist\'orica de esta homolog\'ia puede ser encontrada en \cite{PBIMW}.} \cite{Kho}. Luego el matem\'atico polaco J\'ozef Przytycki utiliz\'o un m\'etodo diferente para calcular esta homolog\'ia; en particular, en su m\'etodo us\'o la conexi\'on que existe entre la homolog\'ia de Khovanov y la homolog\'ia de Hochschild \cite{Prz3}. En esta secci\'on presentamos una nueva forma para obtener este resultado usando directamente la sucesi\'on larga previamente construida.

\begin{definition}
	Un \textbf{enlace toroidal} del tipo (p,q) (tambi\'en llamado $(p,q)$-enlace), es un enlace equivalente a una curva contenida en un toro est\'andar $T^{2}$. Esta curva envuelve $p$ veces alrededor de la longitud y $q$ veces alrededor del meridiano. Si $p$ y $q$ son primos relativos, entonces el enlace tiene un s\'olo componente y se denomina un \textbf{nudo toroidal}.
\end{definition}

Sea $v= \textbf{\KPB}$ un cruce fijo.\\

\begin{figure}[H]
	\centering
	\includegraphics[width=0.75\linewidth]{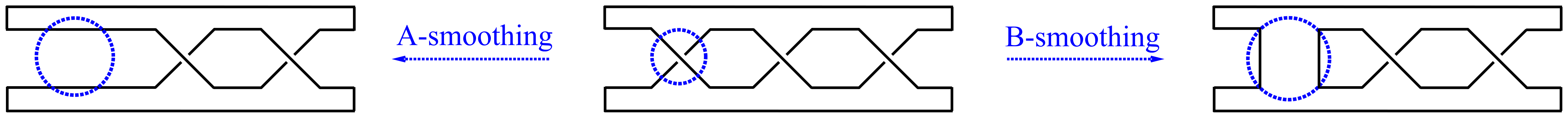}
	\caption{nudo toroidal $T(2,n-1)=T(2,2)$ (izquierda), nudo toroidal $T(2,n)=T(2,3)$ (centro), y el nudo trivial con enmarcado $1-n=-2$ (derecha).}
	\label{fig:torus-23}
\end{figure}

Note que un smoothing vertical en $v$ produce el nudo trivial con $1-n$ ``giros''. Ahora construimos la sucesi\'on larga de la homolog\'ia de Khovanov.
Primero, observe que mantener el cruce (trivialmente) no cambia el diagrama de $T(2,n)$, un smoothing horizontal resulta en el diagrama $T(2,n-1)$, y finalmente, un smoothing vertical resulta en el nudo trivial con ``enmarcado girado'' $1-n$ veces. Por ejemplo, oberve la Figura \ref{fig:torus-23} la cual ilustra este proceso para el nudo toroidal $T(2,3)$.

El siguiente teorema establece la homolog\'ia de Khovanov para nudos toroidales $T(2,n)$ con $n>0$. 

\begin{theorem} 
	Sea $T(2,n)$ un nudo toroidal con $n>0$. Su homolog\'ia de Khovanov $H_{a,b}(T(2,n))$ es dada por:	
	\[ 
	H_{a,b}(T(2,n))=  \left\{   \begin{array}{ll}
		\mathbb{Z} & \mbox{para $(a,b)=(n,n)$ o $(a,b)=(-n,-3n)$},\\
		\mathbb{Z} & \mbox{para $a=n-2s$, $b=n-4s+4$ donde $s$ es par y $0\leq s \leq n$,} \\
		\mathbb{Z} & \mbox{para $a=n-2s$, $b=n-4s$ donde $s$ es impar y $3\leq s \leq n$,} \\
		\mathbb{Z}_2   & \mbox{para $a=n-2s$, $b=n-4s+4$ donde $s$ es impar y $3\leq s \leq n$,}\\
		0           & \mbox{de otra forma.}
	\end{array}
	\right.
	\]	
\end{theorem}

\textbf{Demostraci\'on:} La herramienta principal que usamos en la prueba de este resultado es la sucesi\'on larga de la homolog\'ia de Khovanov (\ref{LES}). La demostraci\'on procede por inducci\'on sobre $n$.	

\ 

\textbf{Para n=1:} El nudo toroidal $T(2,1)$ es el nudo trivial con un ``giro'' positivo. Su homolog\'ia esta dada por:
	\[ H_{a,b}(T(2,1))=
\left\{ \begin{array}{ll}
	\mathbb{Z} & \mbox{para $(a,b)=(1,3\pm 2)$},\\
	0  & \mbox{de lo contrario}.
\end{array}
\right.
\]

Observe que $H_{1,1}(T(2,1))=\mathbb{Z}$ y note que $(a,b)=(1,1)=(n,n)$, lo cual significa que el teorema se cumple. Adem\'as, note que $H_{1,5}(T(2,1))=\mathbb{Z}$ y $(a,b)=(1,5)=(n-2s,n-4s+4)=(1-2(0),1-4(0)+4)$ para $s=0$ y el teorema se cumple. Por lo tanto, conclu\'imos que el teorema se cumple para el caso $n=1$.

\

\textbf{Para n=2:} El nudo toroidal $T(2,2)$ es llamado el enlace de Hopf. El teorema se cumple como puede ser verificado con la Tabla \ref{KHHopf}\footnote{Un resultado bien conocido, vea por ejemplo \cite{MV}.}.	

\begin{table}[h] 
	\begin{minipage}{1\textwidth}
		\centering	
		\begin{tabular}{c||c|c|c|c|c} 
			\ \ $\textbf{b}  \ \ | \ \ \textbf{a}$ \ \  & \ \ \ -2 \ \ \ & \ \ \ 0 \ \ \ & \ \ \ 2 \ \ \  \\
			\hline \hline 
			6                                           &               &               &   $\mathbb{Z}$   \\ 
			\hline
			2                                           &               &                & $\mathbb{Z}$     \\
			\hline
			-2                                           &   $\mathbb{Z}$ &               &              \\
			\hline
			-6                                           & $\mathbb{Z}$  &                 &               \\
			\hline			
			
		\end{tabular}
	\end{minipage}
	\caption{Khovanov homology table of the torus knot $T(2,2)$.}  \label{KHHopf}
\end{table}

Supongamos como \textbf{hip\'otesis inductiva} que el resultado se cumple para $n-1$ donde $n > 2$. Adem\'as considere el caso en el que el mapa $\beta_{\ast}: H_{a,b}(T(2,n)) \longrightarrow H_{a-1,b-1}(T(2,n-1))$ no es necesariamente un isomorfismo. As\'i, analizamos la sucesi\'on larga de la homolog\'ia de Khovanov. Recordamos que, al solucionar un cruce $v$ con un smoothing de tipo $B$, obtenemos el nudo trivial con enmarcado $1-n$, i.e. $T_{B}(2,n)=\bigcirc^{1-n}$. Por lo tanto, su homolog\'ia esta dada por: 

\

\[ H_{x,y}(\bigcirc^{1-n})=
\left\{ \begin{array}{ll}
	\mathbb{Z} & \mbox{para $(x,y)=(1-n,3(1-n)\pm 2)$},\\
	0  & \mbox{de otra forma}.
\end{array}
\right.
\]

\

Esto significa que $\beta_{\ast}$ no es necesariamente un isomorfismo cuando
 $H_{x,y}(\bigcirc^{1-n}) \neq 0$. Esto es, $H_{x,y}(\bigcirc^{1-n}) \neq 0$ cuando $(x,y)=(1-n,3(1-n)-2)=(1-n,1-3n)$ o $(x,y)=(1-n,3(1-n)+2))=(1-n,5-3n).$
 
\ 

\hspace*{0.5cm}\textbf{Caso (i)} Supongamos que $(x,y)=(1-n,1-3n)$:\\

En la sucesi\'on larga en cercan\'ias de $H_{1-n,1-3n}(\bigcirc^{1-n})$ tenemos:\\

\begin{equation*}
	\begin{split}
		\xrightarrow{\makebox[1cm]{}} H_{1-n,-1-3n}(T(2,n-1)) &   \xrightarrow{\makebox[2cm]{$  \textcolor{blue}{\left( \partial^{Conn} \right)}$}} H_{1-n,1-3n}(\bigcirc^{1-n}) \xrightarrow{\makebox[1.8cm]{\textcolor{blue}{$\alpha_{\ast}$}}}\\		 
		H_{-n,-3n}(T(2,n))\xrightarrow{\makebox[1.5cm]{\textcolor{blue}{$\beta_{\ast}$}}} & H_{-n-1,-3n-1}(T(2,n-1)) \xrightarrow{\makebox[1cm]{}}. \\
	\end{split}
\end{equation*}	

Por la hip\'otesis inductiva, $H_{1-n,-1-3n}(T(2,n-1))$ y $H_{-n-1,-3n-1}(T(2,n-1))$ son triviales. Entonces, la sucesi\'on previa se convierte en:\\

\begin{equation*}
	\begin{split}
		\xrightarrow{\makebox[1cm]{}} 0 \xrightarrow{\makebox[2cm]{$  \textcolor{blue}{\left( \partial^{Conn} \right)}$}} H_{1-n,1-3n}(\bigcirc^{1-n}) & \xrightarrow{\makebox[1.8cm]{\textcolor{blue}{$\alpha_{\ast}$}}}\\		 
		H_{-n,-3n}(T(2,n))\xrightarrow{\makebox[1.5cm]{\textcolor{blue}{$\beta_{\ast}$}}} & 0 \xrightarrow{\makebox[1cm]{}}. \\
	\end{split}
\end{equation*}		
	As\'i, $H_{1-n,1-3n}(\bigcirc^{1-n})=H_{-n,-3n}(T(2,n))=\mathbb{Z}$, y por tanto el teorema se cumple pues $(a,b)=(-n,-3n)$.
	
\

\hspace*{0.5cm}\textbf{Caso (ii)} Supongamos que $(x,y)=(1-n,5-3n)$:\\	

En la sucesi\'on larga de homolog\'ia de Khovanov, en cercan\'ias de $H_{1-n,5-3n}(\bigcirc^{1-n})$ tenemos:
\begin{equation*}
	\begin{split}
		0 \xrightarrow{\makebox[1cm]{}} H_{2-n,4-3n}(T(2,n)) \xrightarrow{\makebox[1.5cm]{\textcolor{blue}{$\beta_{\ast}$}}} & H_{1-n,3-3n}(T(2,n-1))	 
		\xrightarrow{\makebox[2cm]{$  \textcolor{blue}{\left( \partial^{Conn} \right)}$}}\\ H_{1-n,5-3n}(\bigcirc^{1-n})\xrightarrow{\makebox[1.5cm]{\textcolor{blue}{$\alpha_{\ast}$}}} H_{-n,4-3n}(T(2,n))& \xrightarrow{\makebox[1cm]{\textcolor{blue}{$\beta_{\ast}$}}} \\
		&H_{-n-1,3-3n}(T(2,n-1)) \xrightarrow{\makebox[1cm]{}}. \\
	\end{split}
\end{equation*}	

Por la hip\'otesis inductiva $H_{-n-1,3-3n}(T(2,n-1))=0$ y $H_{1-n,3-3n}(T(2,n-1))=H_{-(n-1),-3(n-1)}=\mathbb{Z}$. As\'i, la sucesi\'on previa se convierte en:
\begin{equation*}
	\begin{split}
		0 \xrightarrow{\makebox[1cm]{}} H_{2-n,4-3n}(T(2,n)) \xrightarrow{\makebox[1.5cm]{\textcolor{blue}{$\beta_{\ast}$}}}  \mathbb{Z}	 
		\xrightarrow{\makebox[2cm]{$  \textcolor{blue}{\left( \partial^{Conn} \right)}$}} & \\ \mathbb{Z} \xrightarrow{\makebox[1.5cm]{\textcolor{blue}{$\alpha_{\ast}$}}} H_{-n,4-3n}(T(2,n))& \xrightarrow{\makebox[1cm]{\textcolor{blue}{$\beta_{\ast}$}}} 
		0 \xrightarrow{\makebox[1cm]{}}. \\
	\end{split}
\end{equation*}	

Las entradas restantes en la sucesi\'on se determinan entendiento el mapa $\partial^{Conn}_{1-n,3-3n}: \mathbb{Z} \longrightarrow \mathbb{Z}$. \\
En general hay dos posibilidades: el mapa es el mapa cero o es multiplicaci\'on por $k>0$. 

Supongamos primero que $\partial^{Conn}_{1-n,3-3n}$ es el mapa cero. En este caso,
$$H_{2-n,4-3n}(T(2,n))=H_{1-n,3-3n}(T(2,n-1))=\mathbb{Z} \ \ \ \ \ y,$$
$$H_{-n,4-3n}(T(2,n))=H_{1-n,5-3n}(\bigcirc^{1-n})=\mathbb{Z}.$$
As\'i, $H_{2-n,4-3n}(T(2,n))=\mathbb{Z}$ y en este caso tenemos que $(a,b)=(2-n,4-3n)=(n-2n+2,n-4n+4)=(n-2(n-1),n-4(n-1))=(n-2s,n-4s)$ para $s=n-1$. De forma similar, $H_{-n,4-3n}(T(2,n))=\mathbb{Z}$ donde $(a,b)=(-n,4-3n)=(n-2s,n-4s+4)$ para $s=n$ y as\'i, el teorema se cumple.\\

Finalmente, suponga que $\partial^{Conn}_{1-n,3-3n}$ es multiplicaci\'on por $k>0$. Este caso particular fue estudiado en 2006 por M. Pabiniak, J. H. Przytycki, y R. Sazdanović en \cite{PPS} usando el \'algebra $\mathcal{A}_{m}=\quotient{\mathbb{Z}[x]}{(x^{m})}$, la cual para $m=2$ es fuertemente relacionada a la homolog\'ia de Khovanov \cite{Kho, BN}. Cuando $n$ es par, $\partial^{Conn}_{1-n,3-3n}$ es el mapa cero y as\'i $H_{2-n,4-3n}(T(2,n))=\mathbb{Z}=H_{-n,4-3n}(T(2,n))$, como se observ\'o en la parte previa. Cuando $n$ es impar, $\partial^{Conn}_{1-n,3-3n}$ es multiplicaci\'on por $k=2$  y obtenemos que $H_{2-n,4-3n}(T(2,n))=0$, $H_{-n,4-3n}(T(2,n))=\mathbb{Z}_{2}$, y as\'i el teorema se cumple.

\ 

As\'i damos por finalizada la demostraci\'on.

\subsection{Ejemplo: Tablas de homolog\'ia para T(2,11) y T(2,12)}\label{examplesKHtables}
Presentamos en esta secci\'on las tablas de homolog\'ia de Khovanov Tabla \ref{Table 1} y Tabla \ref{Table 2} para $T(2,11)$ y $T(2,12)$, respectivamente. 

\ 

Observe que $H_{a,b}(T(2,n))$ tiene soporte en dos diagonales que contienen $(n, n)$ o $(n,n+4)$. En otras palabras, $H_{a,b}(T(2,n))$ es no trivial solamente para $(a,b)=(n-2s,n-4s)$ o $(a,b)=(n-2s,n-4s +4)$. Adem\'as, $H_{a,b}(T(2,n))$ tiene torsi\'on solamente en grupos para los cuales $(a,b)=(n-2s,n-4s +4)$.

\begin{table}[ht] 
	\begin{minipage}{1\textwidth}
		\centering	
		\begin{tabular}{c||c|c|c|c|c|c|c|c|c|c|c|c|}
			\ \ $\textbf{b} \ \ | \  \ \textbf{a}$  &-11&-9& -7& -5  &  -3  &  -1  &  1 & 3 & 5 & 7& 9 & 11 \\
			\hline \hline
			15                 & & & & & & & & & & & & $\mathbb{Z}$  \\
			\hline
			11                 & & & & & & & & & & & & $\mathbb{Z}$    \\
			\hline
			7                  & & & & & & & & & &     $\mathbb{Z}$ & &  \\
			\hline
			3                  & & & & & & & & &       $\mathbb{Z}_{2}$ & & &  \\
			\hline
			-1                  & & & & & & & &            $\mathbb{Z}$  &  $\mathbb{Z}$ & & & \\
			\hline
			-5                 & & & & & & &  $\mathbb{Z}_{2}$ & & & & & \\
			\hline
			-9                 & & & & & &  $\mathbb{Z}$ & $\mathbb{Z}$ & & & & &    \\
			\hline
			-13                & & & & &   $\mathbb{Z}_{2}$ & & & & & & & \\
			\hline
			-17                & & & &  $\mathbb{Z}$ & $\mathbb{Z}$ & & & & & & &    \\
			\hline
			-21                & & &   $\mathbb{Z}_{2}$ & & & & & & & & & \\
			\hline
			-25                & &  $\mathbb{Z}$ & $\mathbb{Z}$ & & & & & & & & &    \\
			\hline
			-29                &    $\mathbb{Z}_{2}$ & & & & & & & & & & & \\
			\hline
			-33               & $\mathbb{Z}$   &&&&&&&&&&& \\
			\hline
		\end{tabular}
	\end{minipage}
	\caption{Homolog\'ia de Khovanov para el nudo toroidal $T(2,11)$.} \label{Table 1}
\end{table}

\begin{table}[H]
	\begin{minipage}{1\textwidth}
		\centering
		\begin{tabular}{c||c|c|c|c|c|c|c|c|c|c|c|c|c|}
			\ \ $\textbf{b} \ \ | \  \ \textbf{a}$ &-12 &-10&-8& -6& -4  &  -2  &  0  &  2 & 4 & 6 & 8& 10 & 12 \\
			\hline \hline
			16               & & & & & & & & & & & & & $\mathbb{Z}$  \\
			\hline
			12               & & & & & & & & & & & & & $\mathbb{Z}$    \\
			\hline
			8                & & & & & & & & & & &     $\mathbb{Z}$ & &  \\
			\hline
			4                & & & & & & & & & &       $\mathbb{Z}_{2}$ & & &  \\
			\hline
			0               &  & & & & & & & &            $\mathbb{Z}$  &  $\mathbb{Z}$ & & & \\
			\hline
			-4               & & & & & & & &  $\mathbb{Z}_{2}$ & & & & & \\
			\hline
			-8               & & & & & & &  $\mathbb{Z}$ & $\mathbb{Z}$ & & & & &    \\
			\hline
			-12              & & & & & &   $\mathbb{Z}_{2}$ & & & & & & & \\
			\hline
			-16              & & & & &  $\mathbb{Z}$ & $\mathbb{Z}$ & & & & & & &    \\
			\hline
			-20              & & & &   $\mathbb{Z}_{2}$ & & & & & & & & & \\
			\hline
			-24              & & &  $\mathbb{Z}$ & $\mathbb{Z}$ & & & & & & & & &    \\
			\hline
			-28              & &    $\mathbb{Z}_2$ & & & & & & & & & & & \\
			\hline
			-32               & $\mathbb{Z}$   & $\mathbb{Z}$ &&&&&&&&&&& \\
			\hline
			-36              &  $\mathbb{Z}$   &&&&&&&&&&&& \\
			\hline
		\end{tabular}
	\end{minipage}
	\caption{Homolog\'ia de Khovanov para el enlace toroidal $T(2,12)$.}\label{Table 2}
\end{table}	

\section{Conclusi\'on}\label{CONCLU}
La teor\'ia matem\'atica de nudos es un campo muy activo de investigaci\'on. Especialmente desde el descubrimiento (o construcci\'on) del polinomio de Jones, la teor\'ia se revolucion\'o y desarroll\'o r\'apidamente. Adem\'as de brindar una forma sofisticada y elegante de estudiar nudos, este polinomio sac\'o a la luz conexiones con otras ramas de la ciencia como mec\'anica estad\'istica, teor\'ia cu\'antica de campos, \'algebra de operadores, complejidad computacional, biolog\'ia y qu\'imica. Gran parte de la investigaci\'on actual en teor\'ia de nudos est\'a conectada al polinomio de Jones y ahora, adem\'as, incluye el concepto de categorificaci\'on. En particular, las teor\'ias de homolog\'ia en contextos topol\'ogicos generan estructuras matem\'aticas que abren muchas posibilidades para pr\'oximas investigaciones. Por ejemplo, actualmente hay mucha investigaci\'on activa en la homolog\'ia de Khovanov (par e impar), la homolog\'ia de Khovanov-Rozansky, la homolog\'ia de Yang-Baxter, homolog\'ia de quandles, entre otras.

\

Sin duda la mayor parte de la literatura en teor\'ia de nudos est\'a en ingl\'es, aunque existen algunos trabajos en espa\~nol como \cite{MRS, Prz4}. Con miras a atacar esa brecha del idioma, una traducci\'on parcial al ingl\'es de este trabajo es ofrecida en \cite{MV2}. El lector es referido a los trabajos de Bar-Natan y Khovanov para una discusi\'on detallada de la homolog\'ia de Khovanov \cite{BN, Kho}. Para una introducci\'on completa a distintos temas actuales de investigaci\'on en teor\'ia de nudos se recomienda \cite{PBIMW}. En cuanto al uso de la sucesi\'on larga de Khovanov, es pertinente decir que existen futuros posibles caminos de investigaci\'on. Por ejemplo, se puede considerar la familia de nudos ``pretzel'' de tres columnas ya que resolver un cruce de forma vertical, producir\'ia un nudo trivial con cierto n\'umero de ``giros''. Se podr\'ia tambi\'en tratar de aplicar una metodolog\'ia similar para nudos toroidales de tipo $T(3,n)$. Adem\'as, se ha verificado que la sucesi\'on larga de Khovanov puede ser muy \'util en crear nudos cuya homolog\'ia contenga tipos espec\'ificos de torsi\'on \cite{MPSWY}. 

\section*{Acknowledgments}
El autor reconoce orgullosamente el apoyo de parte de la Fundaci\'on Nacional de Ciencia de Estados Unidos NSF (National Science Foundation) a trav\'es del Grant DMS-2212736.

\bibliographystyle{alpha}
\bibliography{bibliografia.bib}

\vspace{1cm}

\noindent Department of Mathematics, The Graduate Center City University of New York, NY, USA, and \newline \noindent Department of Mathematics, University of Puerto Rico at R\'io Piedras, San Juan, PR \\
\noindent Email: \textcolor{blue}{gabrielmontoyavega@gmail.com}

\end{document}